\documentclass{article}
\usepackage{epsfig}
\usepackage{graphics}
\usepackage{color}
\usepackage{amsmath}

\usepackage{amssymb}

\newtheorem{Proposition}{Proposition}
  \newtheorem{Remark}{Remark}
  \newtheorem{Corollary}[Proposition]{Corollary}
  \newtheorem{Lemma}[Proposition]{Lemma}
  
  \newtheorem{Theorem}[Proposition]{Theorem}

\def\e{\mathrm{e}}
\def\i{\mathrm{i}}

\def\Tm{\begin{Theorem}\label}
\def\Pp{\begin{Proposition}\label}
\def\Rm{\begin{Remark}\label}
\def\Lm{\begin{Lemma}\label}
\def\Co{\begin{Corollary}\label}

\def\eTm{\end{Theorem}}
\def\ePp{\end{Proposition}}
\def\eRm{\end{Remark}}
\def\eLm{\end{Lemma}}
\def\eCo{\end{Corollary}}

\def\bfC{{\bf C}}

\makeatletter
\@addtoreset{equation}{section}
\makeatother

\def\Box{{\hfill\hbox{\enspace${\sqre}$}} \smallskip}
\def\sqr#1#2{{\vcenter{\vbox{\hrule height .#2pt
                             \hbox{\vrule width .#2pt height#1pt \kern#1pt
                                   \vrule width .#2pt}
                             \hrule height .#2pt}}}}
\def\sqre{\mathchoice\sqr54\sqr54\sqr{4.1}3\sqr{3.5}3}

\newcommand{\bT}[1]{\begin{theorem}\label{#1}}
\newcommand{\be}[1]{\begin{equation}\label{#1}}
\newcommand{\ba}[1]{\begin{eqnarray}\label{#1}}
\newcommand{\ee}{\end{equation}}
\newcommand{\ea}{\end{eqnarray}}
\newcommand{\bl}[1]{\begin{lemma}\label{#1}}
\newcommand{\bp}[1]{\begin{proposition}\label{#1}}
\newcommand{\br}[1]{\begin{remark}\label{#1}}

\def\bchi{\mbox{\raisebox{.4ex}{\begin{Large}$\chi$\end{Large}}}}

\def\bflam{{\boldsymbol \lambda}}

\def\bfbet{{\boldsymbol \beta}}

\def\erm{\mathrm{e}}
\def\irm{\mathrm{i}}
\def\drm{\mathrm{d}}

\def\RR{\mathbb{R}}

\def\CC{\mathbb{C}}
\def\NN{\mathbb{N}}
\def\ZZ{\mathbb{Z}}
\def\z{\noindent}

\def\bfY{{\bf Y}}

\def\bfl{{\bf l}}

\def\bfV{\mathbf{V}}

\def\lone{{L^1}}

\def\bfa{{\bf a}}
\def\bfB{{\bf B}}
\def\bfA{{\bf A}}

\def\bfy{{\bf y}}
\def\hatby{\tilde{\bfy}}

\def\lap{{\cal L}}

\def\bfk{{\bf k}}

\def\bfm{{\bf m}}

\def\lap{ }
\def\lap{{\cal L}}

\def\bor{{\cal B}}

\def\bfk{{\bf k}}

\def\bfm{{\bf m}}

\def\bflam{{\boldsymbol \lambda}}

\def\bfbet{{\boldsymbol \beta}}

\def\erm{\mathrm{e}}
\def\irm{\mathrm{i}}
\def\drm{\mathrm{d}}
\begin{document}
\hyphenation{trans-series}

\title{On optimal truncation of divergent series solutions
of nonlinear differential systems; Berry smoothing.}

\author{Ovidiu Costin\thanks{Mathematics Department, University of
Chicago, 5734 University Avenue, Chicago IL 60637 and\hfill\break  MSRI, 1000
Centennial Drive, Berkeley CA 94720.  e-mail:
costin\symbol{64}math.uchicago.edu} and Martin D. Kruskal \thanks{
Mathematics Department, Rutgers University, Hill Center--Busch Campus,
New Brunswick, NJ 08903; e-mail: kruskal\symbol{64}math.rutgers.edu}}
\date{ }

\maketitle
\begin{abstract}
  We prove that for divergent series solutions of nonlinear (or
  linear) differential systems near a generic irregular singularity,
  the common prescription of summation to the least term is, if
  properly interpreted, meaningful and correct, and we extend this
  method to transseries solutions.  In every direction in the complex
  plane at the singularity (Stokes directions {\em not} excepted)
  there exists a nonempty set of solutions whose difference from the
  ``optimally'' (i.e., near the least term) truncated asymptotic
  series is of the same (exponentially small) order of magnitude as
  the least term of the series.

  There is a family of generalized Borel summation formulas
  $\mathcal{B}$ which commute with the usual algebraic and analytic
  operations (addition, multiplication, differentiation, etc).  We
  show that there is exactly one of them, $\mathcal{B}_0$, such that
  for any formal series solution $\tilde{f}$, $\mathcal{B}_0(\tilde{f})$
  differs from the optimal truncation of $\tilde{f}$ by at most the
  order of the least term of $\tilde{f}$.

  We show in addition that the Berry (1989) smoothing phenomenon is
  universal within this class of differential systems. Whenever the
  terms ``beyond all orders'' {\em change} in crossing a Stokes line,
  these terms vary smoothly on the Berry scale $\arg(x)\sim
  |x|^{-1/2}$ and the transition is always given by the error
  function; under the same conditions we show that Dingle's rule of
  signs for Stokes transitions holds.

\end{abstract}

\section{Introduction}
\label{sec:Intro}

Summation to the least term (optimal truncation of series) can be
traced back to Cauchy's study of the Gamma function (see Cauchy
(1882))\label{Cauchy}.  The method was applied by Stokes (see Stokes
(1904))\label{Stokes} to solutions of differential equations and was
instrumental in his theory of Stokes' phenomenon.  Optimal truncation
techniques have been greatly improved by Dingle (1973)\label{Dingle A},
and in recent years as a result of the new ideas and methods of
hyperasymptotics introduced by Berry (1989), who also discovered a
surprisingly universal transition in hyperasymptotics, often called
Berry smoothing \label{Berry-smooth2}.  Rigorous results on optimal
truncation and smoothing under various assumptions were proved by
McLeod (1992) (\label{McLeod}second order linear equations), Paris (1992)
\label{Paris1}(linear systems with power coefficients), Olver (1995)
\label{Olver1 A}(Gamma function), Olver and Olde Daalhuis (1995a,b)
\label{Olver-Daalh1} (linear second order systems)
\label{Olver-Daalh2}, Olde Daalhuis (1995, 1996)
\label{Daalhuis} (linear second order systems), the authors (1996)
(\label{CK1 A}first order nonlinear and second order linear) and others.
A wide range of interesting applications were considered by Berry and
Howls
\label{Berry-Howls} \label{Berry-hyp} (1990, 1993)
(saddle point integrals), Berry and Keating (1992) (asymptotics of
the zeta function) and Berry (1991, 1994,
1995)\label{Berry-zet}
\label{Berry-eva}
\label{Berry-gamma}  (asymptotics of the Gamma function, 
applications to dynamical systems, and PDE's).  Other techniques of
exponentially improved asymptotics have found interesting applications
in dynamical systems (see Neishtadt 1984,
\label{Neishtadt} Ramis and Sch\"afke 1996). {\label{Ramis}}

In this paper we prove that asymptotic series of solutions of generic
analytic (linear or nonlinear) differential systems near irregular
singularities are summable to the least term, with errors of the order
of magnitude of the least term with respect to the functions
associated to the series by a specific generalized Borel summation
(the balanced average). In fact we show that the same is true for all
series components of the transseries solutions of such differential
systems.

We show that Berry smoothing extends to this class of asymptotic
series.

\smallskip

We deal with  equations of the form

\begin{eqnarray}
 \label{eqor1}
  \bfy'=\mathbf{f}(x,\bfy)  \qquad \bfy\in\CC^n              
   \end{eqnarray}
for large $x$ in some direction, and study  solutions that 
go to zero in this limit. The assumptions on (\ref{eqor1}) are

(a1) The function $\mathbf{f}$ is analytic
at $(\infty,0)$.

(a2) The eigenvalues $\lambda_i$ of the linearization

\begin{eqnarray}
  \label{linearized}
  \hat{\Lambda}:=-\left(\frac{\partial f_i}{\partial
    y_j}(\infty,0)\right)_{i,j=1,2,\ldots n}
\end{eqnarray}

\z are different from zero and nonresonnant (cf.  \S\ref{nonres}).

It is worth mentioning that a large class of differential systems are
not presented in the form described above but can be brought to it by
suitable changes of variables. For example, under the change of
variables $2x^{3/2}/3=t$, the Airy equation $y''-xy=0$ becomes

$$\frac{\mathrm{d}^2y}{\mathrm{d}t^2}+\frac{1}{3t}\frac{\mathrm{d}y}{\mathrm{d}t}-y=0$$

\z Also, the Painlev\'e I
equation $y''=6y^2+x$ becomes 

$$\frac{\mathrm{d}^2h}{\mathrm{d}t^2}+\frac{1}{t}\frac{\mathrm{d}h}{\mathrm{d}t}-h-\frac{3}{2}h^2-\frac{392}{1875t^4}=0$$

\z by taking $y(x)=-i\sqrt{6x}/2 [1/3-4t^{-2}/75 + h(t)]$ and
$t=(-24x)^{5/4}/30$. In this form, written as systems, both equations
satisfy our assumptions (see \S\ref{nonres} and the notes therein).
The transformations themselves, when they exist, are readily obtained
by comparing the complete formal solution (transseries) of the
original equation with the transseries of the normalized equation,
(\ref{eqformgen,n}) below.

It is  convenient to pull out the inhomogeneous and the linear
terms (relevant to leading order asymptotics) and rewrite the system
in the form

\begin{eqnarray}\label{eqor}
{\bf y}'={\bf f}_0(x)-\hat\Lambda {\bf y}-
\frac{1}{x}\hat B {\bf y}+{\bf g}(x,{\bf y})
\end{eqnarray}

Under the assumptions (a1) and (a2) equation (\ref{eqor}) admits, in normalized
form, an $n-$parameter family of formal exponential series solutions

\begin{eqnarray}
  \label{eqformgen,n}
   \tilde{\bfy}=\tilde{\bfy}_0+\sum_{\bfk\ge 0;|\bfk|>
     0}C_1^{k_1}\cdots C_n^{k_n}
\mathrm{e}^{-(\bfk\cdot\bflam) x}x^{\bfk\cdot\bfm}\tilde{\bfy}_{\bfk}
\end{eqnarray}

\z (see \S\ref{sec:normal}) where
$\tilde{\mathbf{y}}_\bfk=x^{-\bfk(\bfbet+\bfm)}
\sum_{l=0}^{\infty}\bfa_{\bfk;l} x^{-l}$ are formal power series.

We study those formal solutions (\ref{eqformgen,n}) which are at the
same time {\em asymptotic expansions} (or proper transseries, in the
sense of \'Ecalle 1993\label{Ecalle-book A}\label{Ecalle A}). In our context,
only {\em finitely} many powers can be present, of any exponential
which is {\em not small} in the given direction. Note however that any
exponential will decrease in {\em some} direction, and our analysis
eventually covers an $n$ parameter family of formal solutions.

For {\em linear} equations, any formal solutions (\ref{eqformgen,n})
contain finitely many exponentials and are proper transseries.  For
nonlinear equations we require, for large $x$ in a given 
direction $d$ in $\CC$, that $C_i=0$ for any $i$ such that
$\mathrm{e}^{-\lambda_i x}\not \rightarrow 0$ (more details in \S\ref{sec:normal}).

We start by illustrating the various concepts on a very simple
example. The equation $f'+f=1/x$ has the general solution
$y(x)=e^{-x}\mbox{Ei}(x)+C\mathrm{e}^{-x}$
($\mbox{Ei}(x):=PV\int_0^{x}t^{-1}e^t\mathrm{d}t$).  The general {\em formal}
solution for large $x$ is the elementary transseries
$\tilde{y}_C=\sum_{k=0}^{\infty}k!x^{-k-1}+C\erm^{-x}$ with $C\in\CC$
arbitrary. If $\arg(x)\in(\pi/2,3\pi/2)$, the exponential term is
large and classical asymptotics gives full meaning to the family
$\tilde{y}_C$: for each $\tilde{y}_C$ there is a unique actual
solution $y_C$ of the equation such that $y_C\sim\tilde{y}_C$ for large
$x$.  But as $\arg(x)$ crosses an antistokes line $\arg(x)=(k+1/2)\pi$
into a sector where $\erm^{-x}$ decays, $\mathrm{e}^{-x}$ becomes
small {\em beyond all orders} of the divergent series and is undefined
(as part of $\tilde{y}_C$) in classical asymptotics. Classical
asymptotics gives up the exponential and retains the information that
all solutions are asymptotic to $\tilde{y}_0$. In this framework there
is no natural way for asymptotically choosing any privileged solution
as corresponding to $\tilde{y}_0$; $\tilde{y}_0$ becomes more of a
property of the differential equation as a whole, shared by all its
particular solutions.

In terms of
estimating solutions out of $\tilde{y}_0$, classical asymptotics
provides {\em polynomial} precision:
$|y(x)-\tilde{y}_0^{[N]}|\le \mbox{Const.}|x|^{-N-2}$ for large $x$,
where 
 $\tilde{y}_0^{[N]}:=\sum_{k=0}^{N}k!x^{-k-1}$.

 One of the most convenient and frequently used techniques going
 beyond Poincar\'e asymptotics is summation to the least term.  When
 $x$ is very large, the terms of $\tilde{y}_0$ start by decreasing
 rapidly ($k!x^{-k-1}\ll (k-1)!x^{-k}$ if $k\ll x$).  It is tempting
 to keep adding these terms as long as they continue to decrease. The
 least term, for $k=k(x)=\lfloor|x|\rfloor$, evaluates to
 $k(x)!/x^{k(x)+1}\propto x^{-1/2}\erm^{-|x|}$. Even taking into
 account the slight ambiguity as to where to stop adding the terms,
 the numerical precision attained at this stage appears to be high
 enough to measure the small exponential beyond all orders.  It turns
 out that in any direction towards infinity there exists (see e.g.
 Costin and Kruskal 1996\label{CK1 B}) an actual solution of the
 equation which lies within $O(x^{-1/2}\erm^{-|x|})$ of the optimally
 truncated $\tilde{y}_0$ (a nearest solution); this solution, $y_0$,
 is then necessarily unique for any fixed $\arg(x)$, since different
 solutions are separated by $\mbox{Const.}\erm^{-x}$.

 As a manifestation of the Stokes phenomenon, the association between
 $\tilde{y}_0$ and $y_0$ changes with $\xi=\arg x$.  In our example
 $y_0=\mathrm{e}^{-x}\mbox{Ei}(x)\pm \pi\mathrm{i}\mathrm{e}^{-x}$ if
 $\pm\xi>0$ and $y_0=e^{-x}\mbox{Ei}(x)$ when $\xi=0$.

Since the solutions
$y_0$ do not have singularities for large $x$, nothing abrupt should
happen to the terms beyond all orders, either. Michael Berry
discovered that there is indeed an intermediate range $\arg
x\sim|x|^{-1/2}$ where the constant $C$ changes smoothly, in a
surprisingly universal way. 

We prove that Berry smoothing applies to all nonlinear systems
satisfying (a1) and (a2), and the transition function is the same for
all these systems. Technically, we show that Berry smoothing is
governed by the dominant behavior at the singularities in Borel space
nearest to the origin, and these singularities have the same essential
features for all equations within our class. Going beyond
nonresonance, the picture changes. However, properly interpreted, Berry
smoothing might persist in many resonant cases.  See Appendix 2 for a
discussion.

As in the toy model above, in some cases (linear homogeneous second
order and first order equations notably, see Costin and Kruskal
1996\label{CK1 C}) there is {\em only one} solution of the original
differential equation nearest to the divergent series. In general
there is a restricted family of such solutions. Olver's (1964) example
illustrates very well that these special solutions
might not be the obvious ones. At least for this reason it is
important to identify a method of recovering them from their
transseries.

Hyperasymptotics as introduced by Berry leads to a substantial
improvement in precision over optimal truncation but nevertheless ends
up with nonzero errors. Recent very interesting developments of the
original ideas of Berry due to Olde Daalhuis (1996)\label{Olver1 B},
\label{Daalhuis B} suggest that for some linear differential equations
hyperasymptotics may eliminate all errors, but
unfortunately at the price of giving up the convenience of optimal
truncation by adding up vastly more terms beyond the least one in an
extended family of expansions.

At present the only technique which at the same time gives full
account of the complete formal solutions and is widely applicable is
generalized Borel summation (see \'Ecalle 1981 and
1993\label{Ecalle-book B} \label{Ecalle B}, and also Costin 1995 and
1997\label{Costin}\label{Costin++A} for rigorous results in the
context of (\ref{eqor})).  In this paper we show that there is a
natural connection between Borel summation and optimal truncation, and
that optimal truncation is compatible to a unique summation procedure,
the balanced averaging introduced in Costin (1997)\label{Costin++B}.

\subsection{Nonresonance}\label{nonres}
  Let $\theta\in[0,2\pi)$ and
$\tilde{\boldsymbol{\lambda}}=(\lambda_{i_1},...,\lambda_{i_p})$ where
$\left|\arg\lambda_{i_j}-\theta\right|\in(-\pi/2,\pi/2)$ (those
eigenvalues contained in the {\em open} half-plane $H_\theta$ centered along
$\mathrm{e}^{\mathrm{i}\theta}$).  We require that
for {\em any} $\theta$:

(1) $\lambda_{i_1},\ldots,\lambda_{i_p}$ are
$\ZZ$-linearly independent. 

(2) The complex numbers in the set $\{\tilde{\lambda}_{i}-\mathbf{k}
\cdot\tilde{\boldsymbol{\lambda}}\in H_\theta:\bfk\in\NN^p,\,
i=1,...,p\}$ (note: the set is \emph{finite}) have {\em distinct}
directions. These are the Stokes lines $d_{i;\bfk}$.

That the set of $\boldsymbol{\lambda}$ which satisfy (1) and (2) has full
measure follows from the fact that (1) and (2) follow from the condition:

\begin{eqnarray}\label{strongnonr}
  \Big(\bfm,\bfm'\in\ZZ^n ,\ 
\alpha\in\RR\   \mbox{and}\
(\bfm-\alpha\bfm')\cdot\boldsymbol{\lambda}=0
\Big)\ \Rightarrow \Big(\bfm=\alpha\bfm'\Big)
\end{eqnarray}

\z Indeed, if (\ref{strongnonr}) fails, one  of $\Re\lambda_j,
\Im\lambda_j$ is a rational function with rational coefficients of the
other $\Re\lambda_j$ and $\Im\lambda_j$, corresponding to a zero
measure set in $\RR^{2n}$.

{\bf Notes}.  (i) Since only small exponentials are allowed in a
proper transseries, our conditions only restrict those (sets of)
eigenvalues which are associated, through (\ref{eqformgen,n}), to
exponentials that can be simultaneously {\em small} in some
direction. 

(ii) If condition (1) fails, then higher terms in the transseries
(which are present for nonlinear systems) may {\em resonate}, that is
$\bfk\cdot\bflam=\bfk'\cdot\bflam$ for some $\bfk> 0,\bfk'> 0$
($\bfk\ne\bfk'$). This affects the equations characterizing
$\bfy_\bfk$. However, we think that our results should not change in
any substantial way unless in fact 
$\lambda_i=\lambda_j$ for some $i\ne j$. See Appendix 2 for examples
of this latter situation.

(iii) In the Airy and P1 equations discussed in the introduction we
have, after normalization, $\lambda_{1,2}=\pm 1$ and the assumptions
of nonresonnance are satisfied.

(iv) For linear systems, a weaker nonresonnance condition would
suffice for our analysis, namely that $\arg\lambda_i=\arg\lambda_j$
for $i\ne j$ (see Costin 1995). This condition, as well as (2) above,
insures that the Stokes directions are distinct.

\subsection{Normalization and  notation.}\label{sec:normal}

By means of normal form calculations (changes of variables), assuming

(a1) and  (a2), it is possible to prepare (\ref{eqor}) so that
(see Wasow 1968 and Tovbis 1992\label{Wasow A}\label{To1})

(n1)
$\hat\Lambda=\mbox{diag}(\lambda_i)$ and 

(n2) $\hat B=\mbox{diag}(\beta_i)$

\z (To this end, it is {\em not} necessary that the original matrix
$\hat{B}$ in (\ref{eqor}) is diagonalizable).

 For convenience, we rescale $x$ and
reorder the components of $\bfy$ so that:

(n3) $|\lambda_i|\ge 1$ for $i=1...n$, $\lambda_1=1$, and
$\phi_i<\phi_j$ (cf.  (a2)) if $i<j$, where $\phi=\arg(\lambda_i)$.
(To simplify notations, we will formulate some of our results relative
to $\lambda_1$; they can be easily adapted to any other $\lambda_i$,
$|\lambda_i|=1$.)

To unify the treatment, by taking $\bfy=\bfy_1 x^{-N}$ for some $N>0$,
we ensure that

(n4) $\Re(\beta_j)<0,\ j=1,2,\ldots,n$.

\z (There is an asymmetry at this point: the opposite inequality cannot be
achieved, in general, as simply and without violating analyticity at
infinity.)  Finally, through a transformation of the form
$\bfy\leftrightarrow\bfy-\sum_{k=1}^M\bfa_k x^{-k}$ we arrange that

(n5) $ \mathbf{f}_0=O(x^{-M-1})\mbox{ and }\mathbf{g}(x,\bfy)=
O(\bfy^2,x^{-M-1}\bfy) $. We choose $M>1+\max_i\Re(-\beta_i)$.

Under the assumptions (a1) and (a2) and after the preparation (n1)
through (n5), the system (\ref{eqor}) admits an $n$--parameter family
of formal exponential series solutions, (\ref{eqformgen,n}) (see
\label{Iwano} Iwano 1957
\label{Wasow C}and Wasow 1968; in Costin 1997\label{Costin++C}, a brief
derivation in the context of (\ref{eqor}) is given). In
(\ref{eqformgen,n}), $\bfC\in\CC^n$ is an arbitrary vector of
constants, and $m_i=1-\lfloor\beta_i\rfloor$.

The terms of an {\em asymptotic} expansion are 
well-ordered with respect to $\ll$. Thus, (\ref{eqformgen,n}) is
asymptotic in a given direction iff any {\em ascending} chain
$\Re(-\bfk_1\cdot\bflam x)\le \Re(-\bfk_2\cdot\bflam x)\le\ldots$,
$\bfk_i\ne\bfk_j$, in (\ref{eqformgen,n}) is {\em finite} (agreeing to
omit the terms with $C_i=0$). For (\ref{eqor}) the condition reads:

(n6) $\xi+\phi_i\in (-\pi/2,\pi/2)$ for all $i$
such that $C_i\ne 0$. This
 selects the eigenvalues that lie in a half-plane centered at
$\overline{x}$. Without loss of generality, we let $\arg(x)$
vary around $\arg(\lambda_1)=0$.  From now on, ${\bflam}
=(\lambda_{i_1},\ldots,\lambda_{i_{n_1}})$,
$\bfbet=(\beta_{i_1},\beta_{i_2},\ldots,\beta_{i_{n_1}})$,
$\bfm=(m_{i_1},m_{i_2},\ldots,m_{i_{n_1}})$ and
$\bfbet'=\bfbet+\bfm$ where the index is selected
by (n6).

We will henceforth consider that (\ref{eqor}) is presented in prepared
form, and use the designation transseries only for those formal solutions
satisfying (n6).

With respect to Borel summation we are using the results and notation
of Costin (1995 and 1997)\label{Costin B}.  (A summary is presented in
Appendix 1.)

\section{Main results}
\label{sec:mainres}

All the analysis in this paper is based on the structure of
singularities in Borel space of the solutions of (\ref{eqor1}) and
does not otherwise use (\ref{eqor1}).  The results and proofs can be
easily adapted to other functions than solutions of differential
systems, having a similar singularity structure.

To simplify the notation, as mentioned in the introduction, we formulate
all results relative to $\lambda_1=1$; they can be reformulated
without any difficulty relative to any other $\lambda_i$ of modulus
one. We can reduce the discussion to the case when $$S_1\ne 0$$ (cf.
theorem~\ref{AS}; see \S\ref{sec:asder} for a discussion).

\z {\bf (a)} Behavior for large $r\in\NN$ of the coefficients $a_{r,0}$ in
the asymptotic series $\tilde{\bfy}_0$.

Let

\begin{eqnarray}
  \label{defGa}
\Gamma_{r,j}= \frac{\Gamma(r-\beta'_j+1)}
{2\pi\mathrm{i}\,\mathrm{e}^{\mathrm{i}(r+1-\beta_j)\phi_j}} 
\end{eqnarray}

\begin{Proposition}\label{ACO} For $r\rightarrow\infty$ with $r\in\NN$
  we have

\begin{eqnarray}
  \label{asderiv1}
  \mathbf{a}_{r+1,0}=\mathbf{Y}_0^{(r)}(0)=\sum_{j;|\lambda_j|=1}
\left(\mathbf{e}_j S_j+\mathbf{h}_j(r)\right)\Gamma_{r,j}
\end{eqnarray}

\z where $\mathbf{h}_j(r)\sim
r^{-1}\sum_{k=0}^{\infty}\mathbf{h}_{j;k}r^{-k}$ for large $r$. The
minimum of $|\mathbf{a}_{r,0}x^{-r}|$ is reached for $|r-|x||\le
Const.$ where the constant does not depend on $x,r$.

\end{Proposition}

{\sc Remark}. We can write 
$\mathbf{a}_{r,0}=\sum_{j;|\lambda_j|=1}\mathbf{a}_{r,0;[j]}
$ with $\mathbf{a}_{r,0;[j]}\sim \mathbf{e}_j\Gamma_{r,j}$
and thus, defining
$\tilde{\mathbf{y}}_{0;[j]}=\sum_{r}\mathbf{a}_{r,0;[j]}
x^{-r}$ we have

$$\tilde{\mathbf{y}}_0=\sum_{j;|\lambda_j|=1}\tilde{\mathbf{y}}_{0;[j]}$$

\z The terms of each series $\tilde{\mathbf{y}}_{0;[j]}$ have the same
argument (phase) iff $x$ belongs to the $j$-th Stokes line, i.e. iff
$\arg(x)=\phi_j$. This is a precise formulation of {\em Dingle's rule
  of signs} (cf. Dingle 1973\label{Dingle B}, pp.7). For resonant
differential systems the rule needs to be reinterpreted. See also the
Appendix for some further comments.

\smallskip

\z {\bf (b)} Difference between solutions and their optimally truncated
asymptotic series.

\smallskip

\begin{Theorem}\label{MT} 
  
  (i) Let $\bfy$ be a solution of (\ref{eqor}) and
  $x=r\erm^{\irm\xi}$, with $r\in\NN$ and $\xi\in[0,2\pi)$. Then the
  difference between $\bfy$ and the optimal truncation of
  $\tilde{\bfy}_0$ is of the order of magnitude of the least term,

\begin{eqnarray}
  \label{mt1}
  \bfy(x)-\tilde{\bfy}_0^{[r]}(x)=O(\mathbf{a}_{0,r}r^{-r})\ \ \
  (\mbox{for fixed $\xi$, as } r\rightarrow\infty)
\end{eqnarray}

\z if and only if: (1) $\bfy$ is the balanced
Borel sum (see Eq. (\ref{defmedC}) and theorem~\ref{CEQ}) of the
formal solution $\tilde{\bfy}_0$, i.e.,

\begin{eqnarray}
  \label{soleqn2}
  \bfy=\lap\bor_{\frac{1}{2}}\tilde{\bfy}_0+\sum_{|\bfk|> 0}\bfC^{\bfk}\mathrm{e}^{-\bfk\cdot\bflam
    x}x^{-\bfk\cdot\bfbet}\lap\bor_{\frac{1}{2}}\tilde{\bfy}_\bfk
\end{eqnarray}

\z  and (2) ${C}_j= 0$ for all $j$ such that $|\lambda_j|=1$.

(ii) A similar estimate holds for the higher series
in the transseries: For any $\bfk$ and large $r$ we have 

\begin{eqnarray}
  \label{otherser}
 \lap\bor_{\frac{1}{2}}\tilde{\bfy}_\bfk(x)-\tilde{\bfy}_{\bfk}^{[r]}=O(\mathbf{a}_{\bfk,r}r^{-r}) \ \ \ (\mbox{as } r\rightarrow\infty)
\end{eqnarray}

\end{Theorem}

\z {\bf (c)} Berry smoothing. We now focus on the thin parabolic
region (the Berry scale) $|x|\gg 1$, $\arg(x)-\arg(\lambda_j)=\Omega
|x|^{-1/2}$ near a Stokes line, say $\RR^+$ (the Stokes line for
$\lambda_1$).

Theorem~\ref{BS} shows that on the Berry scale near a Stokes line, the
constant $C$ (``beyond all orders'') of a given solution changes, and
the transition is  given, as predicted by Berry, by an error function.

\smallskip

\begin{Theorem}\label{BS} Let $x=r\erm^{\irm \Omega r^{-1/2}}$, with
  $\Omega\in\RR$ fixed and $r\in\NN$; as $r\rightarrow \infty$ we have

\begin{multline}
  \label{mt2,2}
   \lap\bor_{\frac{1}{2}}\tilde{\bfy}_0(x)-\tilde{\bfy}_0^{[r]}=\frac{1}{2}S_1\,\mathrm{erf}\left(\frac{\Omega}{\sqrt{2}}\right)\mathbf{y}_{\mathbf{e}_1}+o(e^{-x}x^{-\beta'})=\cr
\frac{1}{2}S_1\,\mathrm{erf}\left(\frac{\Omega}{\sqrt{2}}\right)e^{-x}x^{-\beta'}\mathbf{e}_1+o(e^{-x}x^{-\beta'})
\end{multline}

\end{Theorem}

We note here that there always exist directions along which the
smoothing is visible (namely the directions $\phi_i$ such that
$|\lambda_i|=1$ in our normalization).  For other Stokes directions,
$O(\mathbf{a}_{0,r}r^{-r})$ may be much larger than the contribution
of erfc.  For the directions of $\lambda_j$ with $|\lambda_j|>1$ the
smallest term of the series is too large for the Berry transition to
be seen. The change in the terms beyond all orders when crossing the
Stokes line is always measurable by Borel summation (cf.
theorem~\ref{RE} (ii) and Costin 1997\label{Costin++F}) but might be
too small for optimal truncation.

\smallskip
\z {\bf (d)} Connection with Borel summation.

In the context of generic nonlinear systems there is a one parameter
family of {\em distinct} generalized Borel summation formulas which
are compatible with all usual operations (addition, multiplication,
differentiation, composition and their inverses, when meaningful) and
which preserve reality and asymptotic inequalities (see Costin
1997)\label{Costin++G}.  In other words, the summation operators
cannot be intrinsically distinguished, in terms of their algebraic
and/or classical asymptotics properties.

Nevertheless, proposition~\ref{bormed} shows that only the generalized Borel
summation with $\alpha=0$, the balanced average, is compatible with
least term truncation. Balanced average is thus the only Borel
summation method that guarantees optimal least term truncation
properties (incidentally, the balanced average is also the most
symmetric and ``natural'' in its family). We note that in our context
it can be shown that the balanced average equals the (apparently) more
complicated median average of Ecalle.

\begin{Proposition}\label{bormed} Let $\bfy=\lap\bor_\alpha\tilde{\bfy}_0$.
 Then

 \begin{eqnarray}
   \label{connecb}
   \bfy(r\mathrm{e}^{\mathrm{i}\xi})-\tilde{\bfy}_{0}^{[r]}(r\mathrm{e}^{\mathrm{i}\xi})=O(\mathbf{a}_{0,r}r^{-r})
 \end{eqnarray}

\z for large $r$ and all $\xi$ if and only if  $\alpha=0$.

\end{Proposition}

{\em Note}: with $\tilde{\bfy}_0$ fixed, $\bfy$ is a different
solution in different Stokes sectors, see theorem~\ref{RE}.

\subsection{Proofs and further results}
\label{sec:proofs}

\z To simplify the exposition, detailed proofs are provided for the
principal series $\tilde{\bfy}_0$. There are only minor adjustments
necessary to cover all other components $\tilde{\bfy}_\bfk$ of the
transseries; we mention them at the end.  The notation is the same as
in Costin (1997)\label{Costin++H} and is explained in \S\ref{sec:anset}. We let
$\bfY=\bfY_0$.

Taking

\begin{equation}
  \label{lape}
  \bfy=\int_{\mathcal{C}} \bfY(p)\mathrm{e}^{-xp}dp
\end{equation}

\z where $\mathcal{C}$ is a contour in $\mathcal{R}$ such that
$\Re(xp)>0$ (or a linear combination of such contours as in the
balanced average), we get by successive integrations by parts

\begin{equation}
  \label{asympexp}
  \bfy(x)=\sum_{k=0}^{r}\bfY^{(k)}(0) x^{-k-1}+
x^{-r-1}\int_{\mathcal{C}}\bfY^{(r+1)}(p)\mathrm{e}^{-px}dp
\end{equation}

\z Taking $x=r \mathrm{e}^{ -\mathrm{i}\varphi}$ we get

\begin{equation}
  \label{estimnth}
  \bfy(x)-\sum_{k=0}^{r}\frac{\bfY^{(k)}(0)}{x^{k+1}}
  =r^{-r-1}
 \mathrm{e}^{-\irm(r+1)\varphi}\int_{\mathcal{C}}\bfY^{(r+1)}(p)
 \mathrm{e}^{-pr
 \mathrm{e}^{ -\mathrm{i}\phi}}\drm p
\end{equation}

\z Equation (\ref{estimnth}) represents the error in summation near the
least term; we will show that for large $r$

\begin{eqnarray}
  \label{aim2} \label{aim1}
\int_{\mathcal{C}}\bfY^{(r+1)}(p)
 \mathrm{e}^{-pr
 \mathrm{e}^{ \mathrm{-i}\phi}}\drm p=
 O\left(\bfY^{(r)}(0)\right)
\end{eqnarray}

\subsection{Asymptotic behavior of $\bfY^{\mathbf(r)}$ for large $r\in\NN$}
\label{sec:asder}

(1) If, exceptionally, all Stokes constants are zero (i.e.  $S_j=0$
(cf. (\ref{defSj})) for $j=1,\ldots,n$), then $\bfY$ is entire and
$|\bfY(p)|\le\exp(\nu|p|)$ for some $\nu$ (theorems~\ref{AS} and
contour,

\begin{eqnarray}
  \label{convcase}
  |\bfY^{(r)}(0)|=\left|
 \frac{r!}{2\pi \mathrm{i}}\oint\frac{\bfY(s)}{s^{r+1}}ds\right|\le
\frac{r!}{2\pi }\frac{\exp(\nu r)}{r^{r+1}}2\pi r \sim \sqrt{2\pi r}\mathrm{e}^{r(\nu-1)}
\end{eqnarray}

\z for large $r$. This implies at once that $\hatby_0$ is {\em
  convergent} in $1/x$. This case is thus trivial: summation to the
least term is nothing else than convergent summation, and the sum
is a true solution of the equation.

\z (2)  By theorem~\ref{AS}
for all $i$ such that $S_i= 0$, $\bfY$ is analytic along
$\RR^+\lambda_i$. We thus only count the $\lambda_i$ with $S_i\ne 0$ and rescale variables and
change notations so that $\min\{|\lambda_i|:S_i\ne 0\}=1$ and
$\beta_0=\Re(\beta_1)=\min\{\Re(\beta_i):|\lambda_i|=1,\ \ S_i\ne 0\}$.

\begin{picture}(0,0)%
\epsfig{file=laplacecont.pstex}%
\end{picture}%
\setlength{\unitlength}{0.00041700in}%
\begingroup\makeatletter\ifx\SetFigFont\undefined%
\gdef\SetFigFont#1#2#3#4#5{%
  \reset@font\fontsize{#1}{#2pt}%
  \fontfamily{#3}\fontseries{#4}\fontshape{#5}%
  \selectfont}%
\fi\endgroup%
\begin{picture}(11839,9372)(720,-8725)
\end{picture}

\centerline{Fig. 1}

\begin{Proposition}\label{asymp1} For large $r\in\NN$ we have

\begin{eqnarray}
  \label{nthderiv}
  \bfY^{(r)}(0)
=\sum_j \frac{S_j(r+m_j)!}{2\pi \mathrm{i}}
\int_{\lambda_j}^{\lambda_j(1+\epsilon)}
\frac{\bfY_{\mathbf{e}_j}(s)}{s^{r+m_j+1}}ds+O(r!(1+\epsilon)^{-r})
\end{eqnarray}

\end{Proposition}

{\em Proof}. Let $m=1+\max_j m_j$. Then $\mathcal{P}^m\bfY(p\rm
e^{i\phi})$ is continuous in $\phi$ in $[0,2\pi]\backslash \{\phi_i,i=1...n\}$
and has at most jump discontinuities at $\{\phi_i,i=1...n\}$
(theorem~\ref{CEQ}).  Writing $2\pi\mathrm{i}(r+m)!\bfY^{(r)}(0)=
\oint \mathcal{P}^m\bfY(s)s^{-r-m-1}ds$, pushing the contour of
integration towards the circle of radius $1+\epsilon$ (Fig.
1), and letting $K=\max_\phi \mathcal{P}^m\bfY((1+\epsilon)\rm e^{i\phi})$ and
$N_1=\# \{i:S_i\ne 0\}$ we have

\begin{eqnarray}
  \label{1+eps;estim}
  &&\left|\oint \mathcal{P}^m\bfY(s)s^{-r-m-1}ds-\sum_{i=1}^{N_1}\int_{C_i}\mathcal{P}^m\bfY(s)s^{-r-m-1}ds\right|
\cr &&\le 2\pi K (1+\epsilon)^{-(r+m+1)}
\end{eqnarray}

\z Using 
theorem~\ref{RE}
we have, with $E_{1,2}=O(r!(1+\epsilon)^{-r})$,

\begin{multline}
  \label{pushc}
  \bfY^{(r)}(0)=\frac{(m+r)!}{2\pi\mathrm{i}}\sum_{j=1}^{N_1}
\int_{\lambda_j}^{\lambda_j(1+\epsilon)}S_j\mathcal{P}^{m-m_j}\bfY_{j\mathbf{e}_j}(s-\lambda_j)s^{-r-m-1}ds+E_1\cr=
\sum_{j=1}^{N_1}\frac{(m_j+r)!}{2\pi\mathrm{i}}
\int_{\lambda_j}^{\lambda_j(1+\epsilon)}S_j\bfY_{j\mathbf{e}_j}(s-\lambda_j)s^{-r-m_j-1}ds+E_2
\end{multline}

\z where the $j$-th term has been  integrated by parts $m-m_j$ times
and we used theorem~\ref{AS} to show that the boundary terms
vanish at $s=\lambda_j$, and are $O(r!(1+\epsilon)^{-r})$
at $\lambda_j(1+\epsilon)$. By theorem~\ref{AS} (ii)
we have for small $p$

\begin{eqnarray}
  \label{casesy}
  \bfY_{j\mathbf{e}_j}(p)= p^{\beta'_j-1} {\bfA}_j(p)
\end{eqnarray}

\z and $\bfA_j(0)=\mathbf{e}_j/\Gamma(\beta'_j)$. Changing variables
to $s=\lambda_j(1+z)$ in (\ref{pushc}) we get, by Laplace's method,

\begin{multline}
  \label{intj}
  \frac{(m_j+r)!}{2\pi\mathrm{i}}
  S_j\lambda_j^{\beta'_j-m_j-r-1}\int_0^{\epsilon}z^{\beta'_j-1}(\mathbf{e}_j+
\mathbf{G}_j(z))
  \mathrm{e}^{-(m_j+r+1)\ln(1+z)}\mathrm{d}z\cr=
  \frac{S_j(r+m_j)!}{2\pi\mathrm{i}(r+m_j+1)^{\beta'_j}
\lambda_j^{r+1-\beta_j}}(\mathbf{e}_j+\mathbf{H}_j(r))
\end{multline}

\z where $\mathbf{G}_j(z)=O(z)$ are analytic and $O(z)$ for small $z$
and thus the $\mathbf{H}_j(r)$ have an asymptotic series in $r^{-1}$,
and $\mathbf{H}_j(r)=O(r^{-1})$.  Combining all these estimates we get

\begin{eqnarray}
  \label{finalasympt}
  \bfY^{(r)}(0)
=\sum_{j=1}^{N_1}\frac{\Gamma(r-\beta_j+1)}
{2\pi\mathrm{i}\,\mathrm{e}^{\mathrm{i}(r+1-\beta_j)\phi_j}}(S_j\mathbf{e}_j+O(r^{-1}))
\end{eqnarray}

\z so that, with $x=r\erm^{-\irm\varphi}$,

\begin{multline}
  \label{leasterm2}
  \frac{\bfY^{(r)}(0)}{x^{r+1}}=r^{-r-1}\erm^{\irm
    (r+1)\varphi}\sum_{j;|\lambda_j|=1}\frac{S_j\Gamma(r-\beta_j+1)}
  {2\pi\mathrm{i}\,\mathrm{e}^{\mathrm{i}(r+1-\beta_j)\phi_j}}(\mathbf{e}_j+o(1/r))
\end{multline}

\z Because of the $\mathbf{e}_j$, the  factors inside the sum cannot
cancel each--other. The rest of proposition~\ref{asymp1} is immediate.

\section{Asymptotics of Laplace Integrals}
\label{sec:ALP}

We distinguish two cases:

\subsection{Regular directions}
\label{sec:simplecase}

If $\arg(x)=\xi$ then, by convention, the direction of the Laplace
transform integral is $\arg p=\varphi=-\xi$. 

If $\varphi\ne \phi_i$ then $\bfY$ is analytic and exponentially
bounded along $\arg p=\varphi$ (theorems~\ref{AS} and \ref{CEQ}). We
may assume $\varphi \in(0,\phi_2)$ (cf. \S\ref{sec:asder}).  Let
$M\in\NN$ be large, subject to the following conditions:

($c_1$) $\varphi\pm\arcsin(M^{-1})\in (0,\phi_2)$.

($c_2$) dist$(\{t\erm^{\irm\varphi}:t\ge
M\},\{\RR\lambda_i\}_{i=1..n})>\alpha>1$

($c_3$) $M\cos(\varphi)>2$.

\z Let $\delta_1\ge 2\erm$  be such that
$M+\delta_1\notin\{\NN|\lambda_i|:i=1,2,\ldots,n\}$ (in particular,
there are no singular points
of $\bfY$
on the circle of radius $M+\delta_1$).  Let
$K_1=\max_{t\le M+\delta_1}|\bfY(te^{\irm\varphi})|$
and $K_2=\max_{|s|=M+\delta_1}|\bfY(s)|$. 

Let $N=M+\lfloor\delta_1\rfloor$. By theorem~\ref{AS}, $\bfY(p)$ is
analytic in $\mathcal{J}_M=\{|p|< M+\delta_1,\arg(p)\ne\phi_i\}$ and
by theorem~\ref{CEQ} (i), $\mathcal{P}^{N|\mathbf{m}|}\bfY$ is
uniformly continuous in the interior of $\mathcal{J}_M$; let
$K_3=\max_{|s|=M+\delta_1;j\le|\mathbf{m}|}|\mathcal{P}^{Nj}\bfY(s)|<
(K_1+K_2)(2N|\mathbf{m}|)^{N|\mathbf{m}|}$. By theorems~\ref{AS} and
\ref{CEQ} (ii) we can choose $\nu$ large enough so that for
$\arg(p)\pm\arcsin(M^{-1})\in (0,\phi_2)$

\begin{eqnarray}\label{cddy}
|\mathbf{Y}(p)|\le\exp(\nu(|p|+1))
\end{eqnarray}

We note that  $M,\nu,K_3,N,\delta_1$ are chosen independently of $r$
and of $x$.

By $(c_2)$ and (\ref{cddy}) for 
$t\ge M$ we have 

\begin{eqnarray}\label{cddydiff}
|\bfY^{(r)}(t\erm^{\irm\varphi})|\le \frac{r!}{2\pi}\exp(\nu(t+2))
\end{eqnarray}

\z For
$|x|>M\nu(M-2)^{-1}$ we have

\begin{eqnarray}
  \label{estimMinfi}
  \left| \int_{M\erm^{\irm\varphi}}^{\infty
      \erm^{\irm\varphi}}\erm^{-px}\bfY^{(r)}(p)\drm p\right|\le\frac{r!}{2\pi}
    \frac{\erm^{2\nu}}{\nu-|x|}e^{M(\nu-|x|)}\le K_4  r!\erm^{-2|x|}
\end{eqnarray}

\z  where $K_4$ depends only on $\nu$ and $M$. By construction,
for $t\le M$ we have $|t\erm^{\irm\varphi}-s|>2\erm$ if
$|s|=M+\delta_1$. We then have (see Fig 1.),

\begin{multline}
  \label{nthderiv2}
  \bfY^{(r)}(t\erm^{\irm\varphi})=\frac{(N|\mathbf{m}|+r)!}{2\pi\irm}
  \oint_{\partial\mathcal{J}_M}
  \frac{\mathcal{P}^{N|\mathbf{m}|}\bfY(s)}{(s-t\erm^{\irm\varphi})^
    {r+N|\mathbf{m}|+1}}\drm s
  \cr=\frac{(N|\mathbf{m}|+r)!}{2\pi\irm}\sum_i\int_{C_i}\frac{\mathcal{P}^{N|\mathbf{m}|}\bfY(s)}{(s-t\erm^{\irm\varphi})^
    {r+N|\mathbf{m}|+1}}\drm s+r! \erm^{-2r}E_2\cr=
  \sum_i f_i\int_{C_i}\frac{\mathcal{P}^{Nm_i}\bfY(s)}{(s-t\erm^{\irm\varphi})^
    {r+Nm_i+1}}\drm s+r! \erm^{-2r}E_3
\end{multline}

\z where $$f_i={(Nm_i+r)!}/({2\pi\irm})$$

\z $E_2$ contains the
contribution of the outer circle, and in the last equality we
integrated by parts $|\mathbf{m}|-m_i$ times and included the boundary
terms in $E_3$. We have $|E_2|+|E_3|\le E_4$ for some $E_4$ depending only on
$K_3, M$ and $|\mathbf{m}|$.
In conclusion, 

\begin{multline}
  \label{conc1}
  \int_0^{\infty\erm^{\irm\varphi}} \mathrm{e}^{-px}\bfY^{(r)}(p)dp=
  \sum_{j=1}^n f_j \int_{0}^{M\erm^{\irm\varphi}}\erm^{-px}\drm p\int_{C_j}
  \frac{\mathcal{P}^{m_j
      N}\bfY(s)}{(s-p\erm^{\irm\varphi})^{Nm_j+r+1}}\drm s\cr+
  r!(e^{-2r}+e^{-2|x|})E_5=\sum_{j=1}^n f_j \int_{C_j}\drm s
  \mathcal{P}^{m_j N}\bfY(s)\int_{0}^{M\erm^{\irm\varphi}}
  \frac{\mathrm{e}^{-px}\drm p }{(s-p\erm^{\irm\varphi})^{Nm_j+r+1}}\cr+
  r!(e^{-2r}+e^{-2|x|})E_5
\end{multline}

\z for large enough $r,x$ where $E_5$ is independent of $r,x$.

We have, with $p=t \mathrm{e}^{\mathrm{i}\varphi}, x=r
\mathrm{e}^{-\mathrm{i}\varphi}$, $R=m_j N+1$:

\begin{multline}
  \label{maininte}
  \int_0^{M \mathrm{e}^{\mathrm{i}\varphi}}\frac{ \mathrm{e}^{-px}}{(s-p)^{r+R}}dp=
\int_0^{\infty  \mathrm{e}^{\mathrm{i}\varphi}}\frac{ \mathrm{e}^{-px}}{(s-p)^{r+R}}dp+O( \mathrm{e}^{-M(r+1)})\cr
= \mathrm{e}^{\mathrm{i}\varphi}\int_0^{\infty}\frac{ \mathrm{e}^{-tr}}{(s-t \mathrm{e}^{\mathrm{i}\varphi})^{r+R}}dt+O( \mathrm{e}^{-M(r+1)})
\end{multline}

\z The behavior for large $r$ of the last integral is obtained by
standard steepest descent:  The function $
f(t)=\mathrm{e}^{-t}(\kappa-t)^{-1}$ has a saddle point at
$t=\kappa-1$ and, for $\kappa\notin[0,\infty]$ the steepest descent
path originating at zero starts in a direction opposite to the pole
$\kappa$ and continues to $+\infty$. For $\arg(s)\ne \varphi$, the
contour in (\ref{maininte}) can thus be deformed to a steepest descent
curve without crossing the pole or the saddle.  The main
contribution to the integral comes then from a region where $t$ is
small. We have the
estimate for (\ref{maininte}):

\begin{multline}
  \label{estmaininte}
  \int_0^{\infty}\frac{ \mathrm{e}^{-tr}}{(s-t
    \mathrm{e}^{\mathrm{i}\varphi})^{r+R}}dt =\int_0^{\infty}
  \exp\left(-(r+R) \ln (s-t{ \mathrm{e}^{
      \mathrm{i}\varphi}})-tr)\right)dt\cr=\frac{1}{s^{r+R}}\int_0^{\infty}
  \exp(-t(r-(r+R) \mathrm{e}^{
      \mathrm{i}\varphi}s^{-1})+O(t^{-2}))dt\\=
[rs^{r+R}
    (\mathrm{e}^{\mathrm{i}\varphi}s^{-1}-1)]^{-1}(1+O(r^{-1}))
\end{multline}

\z As $\mathcal{P}^{Nm_j}\bfY^{\pm}$ are continuous along $C_j$,
we obtain

\begin{eqnarray}
  \label{intsup}
  \left|\int_{\lambda_j(1+\epsilon)}^{(M+\delta_1)\erm^{\irm\phi_j}}
s^{-r-R}\frac{\mathcal{P}^{Nm_j}\bfY^{\pm}(s)}{\erm^{\irm\varphi}
s^{-1}-1}\drm s\right|
\le (1+\epsilon)^{-r-R} E_6
\end{eqnarray}

\z for some  $E_6$ independent of $r$ and $x$. Since 
$\mathcal{P}^{m_j}\bfY^{\pm}$ are continuous for small $\epsilon$ on
the segment between $\lambda_j$ and $\lambda_j(1+\epsilon)$
we have, using (\ref{mainresur}), (\ref{intsup}) and integration by parts, with
$c_j$ being the part of $C_j$ within $\epsilon$ distance
of $\lambda_j$,

\begin{multline}
  \label{subsegment}
   f_j \int_{C_j}\drm s
  \mathcal{P}^{m_j N}\bfY(s)\int_{0}^{M\erm^{\irm\varphi}}
  \frac{\mathrm{e}^{-px}\drm p }{(s-p\erm^{\irm\varphi})^{Nm_j+r+1}}\cr=
f_j \int_{c_j}\drm s
  \mathcal{P}^{m_j N}\bfY(s)\int_{0}^{M\erm^{\irm\varphi}}
  \frac{\mathrm{e}^{-px}\drm p }{(s-p\erm^{\irm\varphi})^{Nm_j+r+1}}
+O((1+\epsilon)^{-r})\cr
=f_j \int_{c_j}
  \mathcal{P}^{m_j}\bfY(s)(-\mathcal{P})^{m_j(N-1)}\int_{0}^{M\erm^{\irm\varphi}}
  \frac{\mathrm{e}^{-px}\drm p }{(s-p\erm^{\irm\varphi})^{Nm_j+r+1}}
+O((1+\epsilon)^{-r})\cr
=e_j \int_{\lambda_j}^{\lambda_j(1+\epsilon)}
  S_j\bfY_{\mathbf{e}_j}(s-\lambda_j)\int_{0}^{M\erm^{\irm\varphi}}
  \frac{\mathrm{e}^{-px}\drm p }{(s-p\erm^{\irm\varphi})^{m_j+r+1}}
+O((1+\epsilon)^{-r})\cr=
e_j \int_{\lambda_j}^{\lambda_j(1+\epsilon)}\mathrm{d}s
  \frac{S_j\bfY_{\mathbf{e}_j}(s-\lambda_j)}{rs^{r+m_j+1}
    (\mathrm{e}^{\mathrm{i}\varphi}s^{-1}-1)}(1+O(r^{-1}))
+O((1+\epsilon)^{-r})
\end{multline}

\z with $$e_j={(m_j+r)!}/({2\pi\irm})$$ The last integral
is very similar to (\ref{pushc}) and is estimated in the 
same way, giving

\begin{multline}
  \label{pushc2}
  e_j \int_{\lambda_j}^{\lambda_j(1+\epsilon)}\mathrm{d}s
  \frac{S_j\bfY_{\mathbf{e}_j}(s-\lambda_j)}{rs^{r+m_j+1}
    (\mathrm{e}^{\mathrm{i}\varphi}s^{-1}-1)}=
\cr=\frac{S_j(r+m_j)!}{2\pi\mathrm{i}r^2
(r+m_j+1)^{\beta'_j-1}
\lambda_j^{r+1-\beta_j}(\erm^{\irm\varphi}\lambda_j^{-1}-1)}(\mathbf{e}_j+O(r^{-1}))
\end{multline}

\z where, in view of (\ref{estimnth}) we need to take $r+1$ instead
of $r$. Noting that $(r+m_j)!r^{-\beta'_j-1}=
\Gamma(r-\beta_j+1)(1+O(r^{-1}))$ and comparing
(\ref{pushc2}) to (\ref{finalasympt}), (\ref{aim2})
is proven.

\Box
\subsection{Stokes directions; balanced averages; Berry smoothing}
\label{sec:crit}

\z This section deals with the special and important case of
exponential asymptotics on and near the Stokes line corresponding to
one of the $\lambda_i$ of largest module ($|\lambda_i|=1$ in our
normalization; without loss of generality we take $i=1$).  Berry's
smoothing formula is proved, and we show that the balanced average is
the only Borel summation process compatible with optimal
truncation. The generalized Borel summations that have good algebraic
and analytic properties are given by (\ref{defmedC}).  Using
theorem~\ref{AS} (i) which shows that for $p\in(0,1)$ we have
$\bfY^+_0(p)=\bfY^-_0(p)$, formula
(\ref{mainresur}), and the estimates in theorem~\ref{CEQ} (i) for
$(\lap\bfY_\bfk)_\gamma$ we see that

\begin{eqnarray}
  \label{represaver1}
  \bfy:=\bfy_\alpha=\mathcal{LB}_\alpha\tilde{\bfy}_0=(1-\alpha)\mathcal{L}\bfY_0^++
  \alpha\lap \bfY_0^- +O(\erm^{-(2-\epsilon)x})
\end{eqnarray}

\z Choosing
$0<\delta<\min_i|\arg(\lambda_1)-\arg(\lambda_i)|$,
relation (\ref{estimnth}) reads, in view of theorem~\ref{CEQ},

\begin{equation}
  \label{estimnthmed}
  \bfy(x)-\sum_{k=0}^{r}\frac{\bfY^{(k)}(0)}{x^{k+1}}
  =r^{-r-1}
 \mathrm{e}^{-\irm(r+1)\varphi}\int_{<\alpha;\infty>}\bfY^{(r+1)}(p)
 \mathrm{e}^{-pr
 \mathrm{e}^{ -\mathrm{i}\varphi}}\drm p+O(\erm^{-(2-\epsilon)x})
\end{equation}

\z where $\int_{<\alpha;a>}$ means

$$\alpha\int_0^{a\mathrm{e}^{-\mathrm{i}\delta}}+(1-\alpha)\int_0^{a\mathrm{e}^{\mathrm{i}\delta}}$$

\z In order to prove theorem~\ref{BS} and proposition~\ref{bormed}, we
need to estimate (\ref{estimnthmed}) near the Stokes line.  Replacing
everywhere $\int_0^{M\erm^{\irm\varphi}}$ by $\int_{<\alpha;M>}$ and
taking $\varphi=0$, the calculations leading to (\ref{conc1}) work
without any other change and we get:

\begin{multline}
  \label{conc1med}
  \int_{<\alpha;\infty>} \mathrm{e}^{-px}\bfY^{(r)}(p)dp=
  \sum_{j=1}^n f_j \int_{<\alpha;M>}\erm^{-px}\drm p\int_{C_j}
  \frac{\mathcal{P}^{m_j
      N}\bfY(s)}{(s-p)^{Nm_j+r+1}}\drm s\cr+
  r!E_5(e^{-2r}+e^{-2|x|})=\sum_{j=1}^n f_j \int_{C_j}\drm s
  \mathcal{P}^{m_j N}\bfY(s)\int_{<\alpha;M>}
  \frac{\mathrm{e}^{-px}\drm p }{(s-p)^{Nm_j+r+1}}\cr+
  r!E_5(e^{-2r}+e^{-2|x|})
\end{multline}

\z For $j\ne 1$, $s-p$ does not vanish for $s\in C_j$,
$|\arg(p)|<\delta$, and thus

\begin{multline}
  \label{conc1med2}
   \sum_{j=2}^n f_j \int_{C_j}\drm s
  \mathcal{P}^{m_j N}\bfY(s)\int_{<\alpha;M>}
  \frac{\mathrm{e}^{-px}\drm p
    }{(s-p)^{Nm_j+r+1}}\cr
=\sum_{j=2}^n f_j \int_{C_j}\drm s
  \mathcal{P}^{m_j N}\bfY(s)\int_0^M
  \frac{\mathrm{e}^{-px}\drm p
    }{(s-p)^{Nm_j+r+1}}
\end{multline}

\z and the estimates leading to (\ref{pushc2}) are valid
if we take $\varphi=0$. We get, for $j\ne 1$,

\begin{multline}
  \label{subsegmentmed}
   f_j \int_{C_j}\drm s
  \mathcal{P}^{m_j N}\bfY(s)\int_{0}^{M}
  \frac{\mathrm{e}^{-px}\drm p
    }{(s-p)^{Nm_j+r+1}}\cr=
e_j \int_{\lambda_j}^{\lambda_j(1+\epsilon)}\mathrm{d}s
  \frac{S_j\bfY_{\mathbf{e}_j}(s-\lambda_j)}{rs^{r+m_j+1}
    (s^{-1}-1)}(1+O(r^{-1}))
+O((1+\epsilon)^{-r})
\end{multline}

\z and 

\begin{multline}
  \label{pushc2med}
  e_j \int_{\lambda_j}^{\lambda_j(1+\epsilon)}
  \frac{S_j\bfY_{\mathbf{e}_j}(s-\lambda_j)}{rs^{r+m_j+1}
    (s^{-1}-1)}(1+O(r^{-1}))=
\cr=\frac{S_j(r+m_j)!}{2\pi\mathrm{i}r(r+m_j+1)^{\beta'_j}
\lambda_j^{r+1-\beta_j}(\lambda_j^{-1}-1)}(\mathbf{e}_j+O(r^{-1}))
\end{multline}

\z We are therefore left with the problem of estimating

\begin{equation}
  \label{conc1medp}
    J:=\int_{C_1}\drm s
  \mathcal{P}^{m_1 N}\bfY(s)\int_{<\alpha;M>}
  \frac{\mathrm{e}^{-px}\drm p
    }{(s-p)^{Nm_1+r+1}}
\end{equation}

\z
We consider the case $\Omega\le 0$, the other case being treated
symmetrically.

\begin{Proposition}\label{estimEi} Let  $\delta\in(0,\pi/2)$,
 $r_1\in\NN$,
$x=r\erm^{\irm\Omega r^{-1/2}}$, $\Omega\le 0$,
$s\ge 1$ and $R=r+r_1$. 

i) For large $r\in\NN$

\begin{multline}
  \label{eqEi}
  g_\alpha=\alpha\int_0^{M\erm^{\irm\delta}}\frac{\erm^{-px}\drm p}{(s-p)^R}
+(1-\alpha)\int_0^{M\erm^{-\irm\delta}}\frac{\erm^{-px}\drm
  p}{(s-p)^R}
\cr=
\frac{\irm r^{-1/2}}{s^{R-1}}\left(\int_0^{\infty}
\exp\left\{-\frac{1}{2}t^2+\left[\Omega_r-\irm\sigma+\frac{1}{\sqrt{r}}
(\sigma\Omega_r+\irm r_1)\right]t\right\}
\drm
t+O(r^{-1})\right)\cr-(1-\alpha)\frac{2\pi\irm}{(R-1)!}x^{R-1}\erm^{-xs}
\end{multline}

\z with 

\begin{equation}\label{defOmega}
\Omega_r=-\irm\sqrt{r}(\erm^{\irm\Omega r^{-1/2}}-1);\ \ 
\sigma=\sqrt{r}(s-1)
\end{equation}

ii) When $\sigma$ is small we have

\begin{eqnarray}
  \label{asymp2}
  g_\alpha=\irm\erm^{\Omega^2/2}\sqrt{\frac{\pi}{2r}}\erm^{-\sqrt{r}\sigma}
\left({\rm erf}\left(\frac{\Omega}{\sqrt{2}}+1\right)E_1
-2(1-\alpha)E_2+\erm^{-\Omega^2/2}O(r^{-1})\right)\cr
\end{eqnarray}

\z where $E_{1,2}=(1+O(r^{-1/2},\sigma))$, ${\rm erf}(x) =
2\pi^{-1/2}\int_0^x\exp(-t^2)\drm t$.
\end{Proposition}

{\em Proof.}

\begin{multline}
  \label{m,infty}
  \int_0^{M\erm^{\pm\irm\delta}}\frac{\erm^{-px}\drm p}{(s-p)^R}=
 \frac{1}{s^{R-1}}\int_0^{\frac{M}{s}\erm^{\pm\irm\delta}}\frac{\erm^{-psx}\drm p}{(1-p)^R}\cr=
\frac{1}{s^{R-1}}\left( \int_0^{\infty\erm^{\irm\delta}}\frac{\erm^{-psx}\drm p}{(1-p)^R}
+O(\erm^{-\frac{xM}{s}\erm^{\pm\irm\delta}})\right)
\end{multline}

\z  Furthermore

\begin{multline}
  \label{rotcont}
  \int_0^{\infty\erm^{\pm\irm\delta}}\frac{\erm^{-px}\drm
    p}{(s-p)^R}+\frac{\pi\irm(1\mp 1)}{(R-1)!}x^{R-1}e^{-xs}=
\int_0^{\irm\infty}\frac{\erm^{-px}\drm
  p}{(s-p)^R}=\irm\int_0^{\infty}\frac{\erm^{-\irm px}\drm p}{(s-\irm
  p)^R}\cr
\end{multline}

\z Taking
$a=1/2-\epsilon$ for small $\epsilon$ we have

\begin{eqnarray}
  \label{rotcont2}
  \int_0^{\infty}\frac{\erm^{-\irm psx}\drm p}{(1-\irm
    p)^R}=\int_0^{r^{-a}}\frac{\erm^{-\irm psx}\drm p}{(1-\irm
    p)^R}+O\left(\exp\left(-\frac{1}{2}r^{2\epsilon}\right)\right)
\end{eqnarray}

\z and 

\begin{multline}
  \label{rotcont3}
\int_0^{r^{-a}}{\erm^{-\irm sxt-R\ln(1-\irm
    t)}\drm t}=
\int_0^{r^{-a}}{\erm^{-\irm sxt+R\irm
    t-Rt^2/2}\drm t}+O(\frac{1}{r})\cr=
\int_0^{\infty}{\erm^{-\irm sxt+R\irm
    t-Rt^2/2}\drm t}+O(\frac{1}{r})=
 r^{-1/2}\left(\int_0^{\infty}\erm^{\irm r^{-1/2}z(R-sx)-\frac{R}{2r}z^2}
\drm
z+O\left(\frac{1}{r}\right)\right)\cr
\end{multline}

\z Part (i) of proposition~\ref{estimEi} follows by combining
(\ref{defOmega}...\ref{rotcont3})). Part (ii) is a straightforward
calculation from (i), using Stirling's formula and Laplace method.

\Box

\z We note that with $\Omega, r_1$ held constant the integral in
(\ref{eqEi}) is bounded for $s\ge 1$ uniformly in $r\in\NN$. Thus, for
$s\ge 1+\epsilon$ we have $g_\alpha\le K_2 e^{-r\epsilon}$ for some
$K_2$ and thus we obtain from (\ref{conc1medp}), proceeding as for
(\ref{subsegment}),

\begin{multline}
  \label{conc1medp3}
   J=\int_{C_1;|s|<1+\epsilon}\drm s
  \mathcal{P}^{m_1 N}\bfY(s)\int_{<\alpha;M>}
  \frac{\mathrm{e}^{-px}\drm p
    }{(s-p)^{Nm_1+r+1}}+E_3\cr=
\frac{e_1}{f_1}\int_1^{1+\epsilon}\drm s S_1\bfY_{\mathbf{e}_1}(s-1)\int_{<\alpha;M>}
  \frac{\mathrm{e}^{-px}\drm p
    }{(s-p)^{m_1+r+1}}+E_4\cr=
\frac{e_1 S_1}{f_1}\int_1^{1+\epsilon}\drm s (s-1)^{\beta'_1-1}(\mathbf{e}_1
+(s-1)h(s))\int_{<\alpha;M>}
  \frac{\mathrm{e}^{-px}\drm p
    }{(s-p)^{m_1+r+1}}+E_4\cr
\end{multline}

\z where $E_{3,4}=O(e^{-r\epsilon})$, and $h(s)$ is smooth on
$[1,1+\epsilon]$.
Taking  $\delta$ small and $r$ correspondingly
large we have, in view of
proposition~\ref{estimEi} (i),

\begin{multline}
  \label{cut2}
  J=\frac{e_1 S_1}{f_1}\int_1^{1+r^{-1+\delta}}\!\!\drm s (s-1)^{\beta'_1-1}(\mathbf{e}_1
+(s-1)h(s))\int_{<\alpha;M>}
  \frac{\mathrm{e}^{-px}\drm p
    }{(s-p)^{m_1+r+1}}+E_5
\end{multline}

\z where $E_5=O(\erm^{-r^{\delta}})$. Using now
proposition~\ref{estimEi} (ii) to evaluate
the inner integral we obtain

\begin{multline}
  \label{cut3}
  f_1 J=\irm \frac{(m_1+r)!}{2\pi\irm} S_1\erm^{\Omega^2/2}\sqrt{\frac{\pi}{2r}}
  \frac{1}{r^{\beta'_1+1}}\left({\rm
    erf}(\Omega/\sqrt{2})-1+
2\alpha\right.\cr\left.+
({\rm
    erf}(\Omega/\sqrt{2})+1) E_6
-2(1-\alpha)E_7+\erm^{-\Omega^2/2}O(r^{-1})\right)
\end{multline}

\z where $E_{6,7}=O(r^{-1/2})$, where, as in (\ref{pushc2}) we need to
take $r+1$ instead of $r$ (cf. (\ref{aim2})).  Taking $\lambda=1$ in
the expression (\ref{pushc2}), which was shown to be of the order of
the least term, we see that the least term is $O(r^{-1/2})$ \emph{smaller}
than (\ref{cut3}) unless

\begin{eqnarray}
  \label{conclusion1}
  {\rm
    erf}(\Omega/\sqrt{2})-1+
2\alpha=0
\end{eqnarray}

(a) Asymptotics along straight lines:  If
$\Omega=0$ ($\arg(x)=0$) we see that only 

$$\alpha=1/2$$

\z (corresponding to the balanced average, which
on the interval $(0,2)$ is the half sum
of the upper and lower continuations) ensures
errors in optimal truncation of the order 
of the least term.

(b) Asymptotics along parabolas and Berry smoothing. Theorem~\ref{BS}
is straightforward application of Stirling's formula to (\ref{cut3}).
(Note also that by theorem~\ref{CEQ} (ii) and theorem~\ref{AS} (ii)
$\mathbf{y}_{\mathbf{e}_1}\sim x^{-\beta'}e^{-x}\mathbf{e}_1$.)
Convergence of the Puiseux series near the singularities of $\bfY_0$
in Borel space was the key in proving universality of the Berry
transition, since using this convergence, the calculation reduced in
effect to the case of one particular function, $\mathrm{Ei}(x)$.

\centerline{*}

The same line of proof as for  $\tilde{\bfy}_0$ works
for $\tilde{\bfy}_\bfk$  as
well. Indeed, we note that
$x^{\bfk\cdot(\bfbet+\bfm)}\tilde{\bfy}_\bfk$ have the same
singularity structure in Borel space as $\tilde{\bfy}_0$. This follows
from theorem~\ref{AS} (ii) and the fact that convolution with a
locally analytic function preserves the convergence of Puiseux series
(cf. Example (3a) in  appendix 1). \Box

\section{Appendix 1.}
\label{sec:app}

\subsection{Results on Borel summation}
\label{sec:anset}

In this paper we make use of some of the results in Costin
(1997)\label{Costin++I}; for convenience we summarize them below.

\z We use the convention $\NN=\NN\cup\{0\}$.  Let 

\begin{eqnarray}
  \label{defW}
  \mathcal{W}=\left\{p\in\CC:p\ne k\lambda_i\,,\forall
k\in\NN,i=1,2,\ldots,n\right\}
\end{eqnarray}

The directions $d_j=\{p:\arg(p)=\phi_j\}, j=1,2,\ldots,n$ are the {\em
  Stokes lines} (Note: sometimes known as {\em anti-}Stokes lines!).  We
construct over $\mathcal{W}$ a surface $\mathcal{R}$, consisting of
homotopy classes of smooth curves in $\mathcal{W}$ starting at the
origin, moving away from it, and crossing at most one Stokes line, at
most once (Fig. 2):

\begin{eqnarray}\label{defpaths}
{\cal R}:=\Big\{\gamma:(0,1)\mapsto \mathcal{W}:\ 
\gamma(0_+)=0;\ \frac{\mathrm{d}}{\mathrm{d}t}|\gamma(t)|>0;\ \arg(\gamma(t))\ \mbox{monotonic}\Big\}\cr
\end{eqnarray} 

\z The Laplace transform along a direction $\phi$ of a function
$F$ 
$\mathcal{L}_\phi F$ will depend in general on $\phi$; the usual
convention is to choose $\phi$ so that $xp\in\RR^+$.  Thus, the Borel
sum of $\tilde{f}$ in the direction $x$, if it exists, is defined as
$\lap_{\phi(x)}\mathcal{B}f$ with $\phi(x)=-\arg(x)$.
By (n6) and the agreed association between $p$ and $x$, and since
Laplace integrals will not depend on the direction of $p$ until
one of the Stokes lines is crossed, we may as well assume that the
direction of integration is either $d_j$ or is arbitrarily close to it.
Define $\mathcal{R}_1$ as the restriction of $\mathcal{R}$ to
$\arg(\gamma)\in(\psi_n-2\pi,\psi_2)$ where $\psi_n=\max\{-\pi/2,\phi_n-2\pi\}$
and $\psi_2=\min\{\pi/2,\phi_2\}$. 

\begin{picture}(0,0)%
\epsfig{file=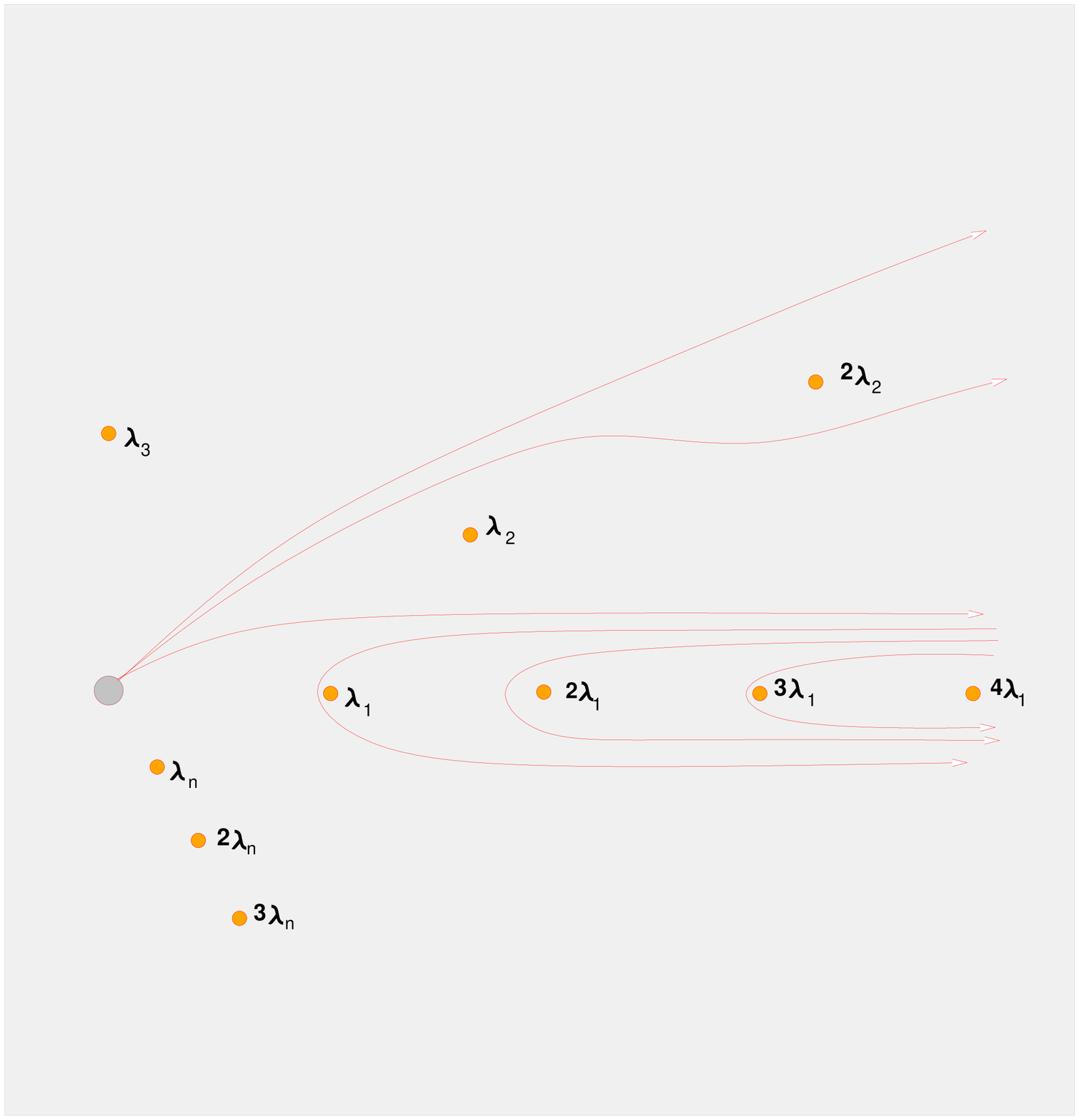, height=6.5cm}%
\end{picture}%
\setlength{\unitlength}{0.00033300in}%
\begingroup\makeatletter\ifx\SetFigFont\undefined%
\gdef\SetFigFont#1#2#3#4#5{%
  \reset@font\fontsize{#1}{#2pt}%
  \fontfamily{#3}\fontseries{#4}\fontshape{#5}%
  \selectfont}%
\fi\endgroup%
\begin{picture}(10890,7989)(4801,-9310)
\end{picture}

\nobreak \smallskip
\centerline{{{Fig 2.} \emph{The paths
near $\lambda_2$ belong to $\mathcal{R}$.   }}}\nobreak
\centerline{\em The paths
near $\lambda_1$ relate to the balanced average}

\bigskip

 We denote the analytic
continuation of $f$ along a curve $\gamma$ by $AC_\gamma(f)$.  For the
analytic continuations near a Stokes line $d_i$ we use notations
similar to \'Ecalle's: $f^-$ is the branch of $f$ along a path $\gamma$
with $\arg(\gamma)<\phi_i$, while $f^{-k+}$  denotes the branch
along a path that crosses the Stokes line between $k\lambda_i$ and
$(k+1)\lambda_i$. We use the notations $\mathcal{P}f$ 
for $\int_0^p f(s)\mathrm{d}s$ and $\mathcal{P}_\gamma f$ if integration
is along a curve $\gamma$. 

We write $\bfk\ge\bfk'$ if $k_i\ge k'_i$ for all $i$ and
$\bfk\succ\bfk'$ if $\bfk\ge\bfk'$ and $\bfk\ne\bfk'$. The relation
$\succ$ is a well ordering on $\NN^{n_1}$.  We let $\mathbf{e}_j$ be
the unit vector in the $j^{\rm th}$ direction.

Formal expansions are denoted with a tilde and capital letters
$\bfY,\bfV\ldots$ will usually denote Borel transforms or functions
otherwise associated to Borel space. For notational convenience, we
will not however distinguish between
$\tilde{\bfY}_k=\bor\tilde{\bfy}_\bfk$, which turn out to be
convergent series, and the sums of these series $\bfY_\bfk$ as germs
of ramified analytic functions.

We have

\begin{eqnarray}\label{Taylor series}
{\bf g}(x,{\bfy})=\sum_{|{\bf l}|\ge 1}{\bf g}_{\bf l}(x) {\bf
y}^{\bf l}=\sum_{s\ge 0;|{\bf l}|\ge 1}{\bf g}_{s,\bf l}x^{-s}
{\bfy}^{\bf l} \ \ (|x|>x_0,|\bfy|<y_0)
\end{eqnarray}

\z where ${\bfy}^{\bf l}=y_1^{l_1}\cdots y_n^{l_n}$ and 
$|{\bf l}|=l_1+\cdots+l_n$. By construction ${\bf g}_{s,\bfl}=0$ 
if $|\bfl|=1$ and $s\le M$.

         \begin{Theorem}\label{AS} (i) $\bfY_0=\bor\tilde{\bfy}_0$
is analytic in $\mathcal{R}\cup \{0\}$. The singularities of $\bfY_0$
(which are contained in the set
$\{l\lambda_j:l\in\NN^+,j=1,2,\ldots,n\}$) are described as follows.
For $l\in\NN^+$ and small $z$

 \begin{multline}
   \label{SY0}
   \bfY_0^{\pm}(z+l\lambda_j)=\pm\Big[(\pm S_j)^l\ln(z)^{0,1}
   \bfY_{l\mathbf{e}_j}(z)\Big]^{(lm_j)}+\bfB_{lj}(z)=\cr
   \Big[z^{l\beta_j'-1}\ln
   z^{0,1}\,\bfA_{lj}(z)\Big]^{(lm_j)}+\bfB_{lj}(z) \ (l=1,2,\ldots)
 \end{multline}

 \z where the power of $\ln(z)$ is one iff
 $l\beta_j\in\ZZ$, and $\bfA_{lj},\bfB_{lj}$ are analytic for small
 $z$. The functions $\bfY_\bfk$ are, in addition, analytic at
 $p=l\lambda_j$, $l\in\NN^+$, iff, exceptionally,

\begin{eqnarray}\label{defSj}
S_j=r_j\Gamma(\beta'_j)\left(\bfA_{1,j}\right)_j(0)=0
\end{eqnarray}

\z where $r_j=1-\mathrm{e}^{2\pi \mathrm{i}(\beta'_j-1)}$
if $l\beta_j\notin\ZZ$ and $r_j=-2\pi \mathrm{i}$
otherwise. The $S_j$ are Stokes constants,
see theorem~\ref{Stokestr}.

(ii) $\bfY_\bfk=\bor{\tilde{\bfy}}_\bfk$, $|\bfk|>1$, are analytic in
$\mathcal{R}\backslash
\{-\bfk'\cdot{\bflam}+\lambda_i:\bfk'\le\bfk,1\le i\le n\}$.  For
$l\in\NN$ and $p$ 
 near $l\lambda_j$, $j=1,2,\ldots,n$ there exist $\bfA=\bfA_{\bfk jl}$ and
$\bfB=\bfB_{\bfk jl}$ analytic at zero so that ($z$ is as above)

 \begin{multline}
   \label{SYK}
  \bfY_\bfk^{\pm}(z+l\lambda_j)=
\pm\Big[(\pm S_j)^l{k_j+l\choose
     l}\ln(z)^{0,1}
   \bfY_{\bfk+l\mathbf{e}_j}(z)\Big]^{(lm_j)}+l\bfB_{\bfk lj}(z)=\cr
  \Big[z^{\bfk\cdot{\bfbet}'+l\beta_j'-1}(\ln z)^{0,1}\,\bfA_{\bfk l
    j}(z)\Big]^{(lm_j)}+l\bfB_{\bfk l j}(z)\ (l=0,1,2,\ldots)
       \end{multline}

       \z where the power of $\ln z$ is $0$ iff $l=0$ or
       $\bfk\cdot{\bfbet}+l\beta_j-1\notin\ZZ$ and $\bfA_{\bfk 0
         j}=\mathbf{e}_j/\Gamma(\beta'_j)$.  Near
       $p\in\{-\bfk'\cdot{\bflam}:0\prec\bfk'\le\bfk\}$, (where
       $\bfY_0$ is analytic) $\bfY_\bfk,\,\bfk\ne 0$ have convergent
       Puiseux series.

\end{Theorem} 

 Let $\bor\tilde{\bfy}_\bfk$ be extended along $d_j$ by
the ``balanced average'' of analytic continuations

\begin{eqnarray}
\label{defmed}
\bor\tilde{\bfy}_\bfk=\bfY_\bfk^{ba}=
\bfY_\bfk^++\sum_{j=1}^{\infty}\frac{1}{2^j}\left(\bfY_\bfk^{-}
-\bfY_\bfk^{-({j-1})+}\right)
\end{eqnarray}

\z The sum above coincides with the one in which $+$ is exchanged with
$-$, accounting for the
reality-preserving property. Clearly, if $\bfY_\bfk$ is analytic along
$d_j$, then the terms in the infinite sum vanish and
$\bfY_\bfk^{ba}=\bfY_\bfk$; we also let $\bfY_\bfk^{ba}=\bfY_\bfk$ if
$d\ne d_j$, where again $\bfY_k$ is analytic. It follows from
(\ref{defmed}) and theorem~\ref{CEQ} below that the Laplace integral
of $\mathbf{Y}^{ba}_\bfk$ along $\RR^+$ can deformed into contours
as those depicted in Fig. 2, with weight $2^{-k}$ for a contour turning around
$(k+1)\lambda_1$. More generally, we consider the averages

\begin{eqnarray}
\label{defmedC0}
\bor_\alpha\tilde{\bfy}_\bfk=\bfY_\bfk^{\alpha}=
\bfY_\bfk^++\sum_{j=1}^{\infty}\alpha^j\left(\bfY_\bfk^{-}
-\bfY_\bfk^{-({j-1})+}\right)
\end{eqnarray}

\z and correspondingly

\begin{eqnarray}
  \label{defmedC}
(\mathcal{LB})_\alpha\tilde{\bfy}_\bfk:=\mathcal{L}\bfY_\bfk^{\alpha}
\end{eqnarray}

\z With $\alpha\in\RR$, this represents the most general family of
averages of Borel summation formulas which commute with complex
conjugation, with the algebraic and analytic operations and have good
continuity properties (see Costin 1995\label{Costin C}).  The value
$\alpha=1/2$ is special in that it is the only one compatible with
optimal truncation.

\begin{Theorem}\label{CEQ} (i) The branches of $(\bfY_\bfk)_\gamma$ in $\mathcal{R}_1$
  have limits in a $C^*$-algebra of distributions,
  $\mathcal{D}'_{m,\nu}(\RR^+)\subset\mathcal{D}'$ (cf.
  \S~\ref{sec:NC}) Their Laplace transforms in
  $\mathcal{D}'_{m,\nu}(\RR^+)$ $\lap(\bfY_\bfk)_\gamma$ exist
  simultaneously and with $x\in\mathcal{S}_x$ and for any $\delta>0$
  there is a constant $K$ and an $x_1$ large enough, so that for
  $\Re(x)>x_1$ we have $\left|\lap(\bfY_\bfk\right)_\gamma(x)|\le
  K\delta^{|\bfk|}$.

In addition, $\mathbf{Y}_\bfk(p\mathrm{e}^{\mathrm{i}\phi})$ are continuous
in $\phi$ with respect to the $\mathcal{D}'_{m,\nu}$ topology,
(separately) on $[\psi_n-2\pi,0]$ and $[0,\psi_2]$.

If
$m>\max_i(m_i)$ and $l<\min_i |\lambda_i|$ then
$\mathbf{Y}_0(p\mathrm{e}^{\mathrm{i}\phi})$ is continuous in
$\phi\in[0,2\pi]\backslash\{\phi_i:i\le n\}$ in the
$\mathcal{D}'_{m,\nu}(\RR^+,l)$ topology and has (at most) jump
discontinuities for $\phi=\phi_i$. For each $\bfk$, $|\bfk|\ge 1$ and
any $K$ there is an $l>0$ and an $m$ such that
$\mathbf{Y}_k(p\mathrm{e}^{\mathrm{i}\phi})$ are continuous in
$\phi\in[0,2\pi]\backslash\{\phi_i; -\bfk'\cdot\bflam+\lambda_i:i\le n
,\bfk'\le\bfk\}$ in the $\mathcal{D}'_{m,\nu}((0,K),l)$ topology and
have (at most) jump discontinuities on the boundary.

(ii) The sum (\ref{defmed}) converges in $\mathcal{D}'_{m,\nu}$ (and
coincides with the analytic continuation of $\bfY_\bfk$ when
$\bfY_\bfk$ is analytic along $\RR^+$). For any
$\delta$ there is a large enough $x_1$ {\em
  independent of $\bfk$} so that  $\bfY^{ba}_\bfk(p)$ with
$p\in\mathcal{R}_1$ are
Laplace transformable
for
$\Re(xp)>x_1$ and furthermore $|(\lap\bfY^{ba}_\bfk)(x)|\le
\delta^{|\bfk|}$. In addition, if $d\ne\RR^+$, then for large $\nu$,
$\bfY_\bfk\in L^1_\nu(d)$.

The functions
$\lap\bfY_\bfk^{ba}$ are analytic for $\Re(xp)>x_1$. For any
$\bfC\in\CC^{n_1}$ there is an $x_1(\bfC)$ large enough so that the
sum

\begin{eqnarray}
  \label{soleqn}
  \bfy=\lap\bfY_0^{ba}+\sum_{|\bfk|> 0}\bfC^{\bfk}\mathrm{e}^{-\bfk\cdot\bflam
    x}x^{-\bfk\cdot\bfbet}\lap\bfY_\bfk^{ba}
\end{eqnarray}

\z converges uniformly for $\Re(xp)>x_1(\bfC)$, and $\bfy$ is a solution
of (\ref{eqor}).  When the direction
of $p$ is not the real axis then, by definition, $\bfY^{ba}_\bfk=\bfY_\bfk$,
$\mathcal{L}$ is the usual Laplace transform 
and (\ref{soleqn}) becomes

\begin{eqnarray}
  \label{soleqnpm}
  \bfy=\lap\bfY_0+\sum_{|\bfk|> 0}\bfC^{\bfk}\mathrm{e}^{-\bfk\cdot\bflam
    x}x^{-\bfk\cdot\bfbet}\lap\bfY_\bfk
\end{eqnarray}

In addition, $\lap\bfY_\bfk^{ba}\sim \tilde{\bfy}_\bfk$ for large $x$
in the half plane $\Re(xp)>x_1$, for all $\bfk$, uniformly.

iii) The general
solution of (\ref{eqor}) that is asymptotic to $\tilde{\bfy}_0$ for
large $x$ along a ray
in $S_x$ can be equivalently written in  the form
(\ref{soleqn}) or as
\begin{eqnarray}
  \label{soleqnp}
  \bfy=\lap\bfY_0^{\pm}+\sum_{|\bfk|> 0}\bfC^{\bfk}\mathrm{e}^{-\bfk\cdot\bflam
    x}x^{-\bfk\cdot\bfbet}\lap\bfY_\bfk^{\pm}
\end{eqnarray}

\z for some $\bfC$ (depending on the solution
and chosen form). With the convention
binding the directions of $x$ and $p$ and the representation form
being fixed, (cf. the beginning of \S\ref{sec:anset}))
the representation of a solution is unique.

\end{Theorem}

\begin{Theorem}\label{RE}
i) For all $\bfk$ and $\Re(p)>j,\Im(p)>0$ as well
as in $\mathcal{D}'_{m,\nu}$ we have

\begin{eqnarray}
  \label{mainresur}
  \bfY_{\bfk}^{\pm j\mp}(p)-\bfY_{\bfk}^{\pm (j-1) \mp}(p) = (\pm S_1)^j\binom{k_1+j}{j}
  \left(\bfY^\pm_{\bfk+j\mathbf{e}_1}(p-j)\right)^{(mj)}
\end{eqnarray}

\z and also, 

\begin{eqnarray}
  \label{thirdresu}
  \mathbf{Y}_{\bfk}^\pm=\bfY_\bfk^{\mp}+\sum_{j\ge 1} {j+k\choose k}(\pm S_1)^{j}(\mathbf{Y}^\mp_{\bfk+j\mathbf{e}_1}(p-j))^{(mj)}
\end{eqnarray}

ii) {\em Local Stokes transition.}

\z Consider the expression of a fixed solution $\bfy$ of (\ref{eqor})
as a Borel summed transseries (\ref{soleqn}). As $\arg(x)$ varies,
(\ref{soleqn}) changes only through $\bfC$, and that change occurs
when the Stokes lines are crossed. We have, in the neighborhood
of $\RR^+$, with $S_1$
defined in (\ref{defSj}):

\begin{eqnarray}
  \label{microsto}
  \bfC(\xi)=\left\{\begin{array}{lll} \bfC^-=\bfC(-0)\qquad&\mbox{for
    $\xi<0$}\\ \bfC^0=\bfC(-0)+\frac{1}{2}S_1\mathbf{e}_1\qquad&\mbox{for
    $\xi=0$}\\ \bfC^+=\bfC(-0)+S_1\mathbf{e}_1\qquad&\mbox{for
    $\xi>0$}\end{array}\right.
\end{eqnarray}
\end{Theorem}

\begin{Remark}\label{R1} In view of (\ref{mainresur}) the different analytic
continuations of $\bfY_0$ along paths crossing
$\RR^+$ at most once can be expressed in terms of
$\bfY_{j\mathbf{e}_1}$. The most general formal solution of (\ref{eqor})
that can be formed in terms of $\bfY_{j\mathbf{e}_j}$ with $j\ge 0$ is
(\ref{eqformgen,n}) with $C_1=\alpha$ arbitrary and $C_j=0$ for $j\ne
1$. Any true solution of (\ref{eqor}) based on such a transseries
is given in (\ref{soleqnp}) with 
$\bfC$ as above. Any average $\mathcal{A}\bfY_0$  along paths
going forward in $\RR^+$ such that 
$\lap \mathcal{A}\bfY_0$ is thus of the form (\ref{defmedC}).
\end{Remark}

\begin{Theorem}\label{Stokestr} Assume only $\lambda_1$ lies
  in the right half plane. Let $\gamma^{\pm}$ be two paths in the
  right half plane, near the positive/ negative imaginary axis such
  that $|x^{-\beta_1+1}\mathrm{e}^{-x\lambda_1}|\rightarrow 1$ as
  $x\rightarrow\infty$ along $\gamma^{\pm}$. Consider the solution
  $\bfy$ of (\ref{eqor}) given in (\ref{soleqn}) with
  $\mathbf{C}=C\mathbf{e}_1$ and where the path of integration is
  $p\in\RR^+$. Then

\begin{eqnarray}\label{classicS}
\bfy=
(C\pm\frac{1}{2}S_1)\mathbf{e}_1 x^{-\beta_1+1}\mathrm{e}^{-x\lambda_1}(1+o(1))
\end{eqnarray}

\z for large $x$ along $\gamma^{\pm}$, where $S_1$ is the same as in
(\ref{defSj}), (\ref{microsto}).

\end{Theorem}

\begin{Proposition}\label{asymptrick} i)  Let $\bfy_1$ and $\bfy_2$ be solutions of (\ref{eqor})
so that $\bfy_{1,2}\sim\tilde{\bfy}_0$ for large $x$ in an
open sector $S$ (or in some direction
$d$);
then $\bfy_1-\bfy_2=\sum_{j}C_j\mathrm{e}^{-\lambda_{i_j}
  x}x^{-\beta_{i_j}}(\mathbf{e}_{i_j}+o(1))$ for some constants $C_j$, where the indices run over
  the eigenvalues $\lambda_{i_j}$ with the property $\Re(\lambda_{i_j}
  x)>0$ in $S$ (or $d$). 
If $\bfy_1-\bfy_2=o(\mathrm{e}^{-\lambda_{i_j}
    x}x^{-\beta_{i_j}})$ for all $j$, then $\bfy_1=\bfy_2$.
 
ii) Let $\bfy_1$ and $\bfy_2$ be solutions of
(\ref{eqor1}) and assume that $\bfy_1-\bfy_2$ has differentiable
asymptotics of the form 
$\mathbf{K}a\exp(-ax)x^b(1+o(1))$ with  $\Re(ax)>0$ and $\mathbf{K}\ne 0$, for large $x$.
Then $a=\lambda_i$ for some $i$.

iii) Let $\mathbf{U}_\bfk\in\mathcal{T}_{\{\cdot\}}$ for all $\bfk$,
$|\bfk|>1$.  Assume in addition that for
large $\nu$ there is a function $\delta(\nu)$ vanishing as
$\nu\rightarrow\infty$ such that

\begin{gather}
\label{cd1}
\sup_{\bfk}\delta^{-|\bfk|}\int_{d}\left|\mathbf{U}_{\bfk}(p)\mathrm{e}^{-\nu p}\right|\mathrm{d}|p|<K<\infty
\end{gather}

\z Then, if $\bfy_1,\bfy_2$ are solutions of (\ref{eqor}) in $S$ where in addition

\begin{eqnarray}
  \label{lapcond}
  \bfy_1-\bfy_2=\sum_{|\bfk|>1}
\mathrm{e}^{-\bflam\cdot\bfk x}x^{\bfm\cdot\bfk}\int_d\mathbf{U}_\bfk(p)
\exp(-xp)\mathrm{d}p
\end{eqnarray}

\z where $\bflam,x$ are as in (n6), then $\bfy_1=\bfy_2$, and
$\mathbf{U}_\bfk=0$ for all $\bfk$, $|\bfk|>1$.

\end{Proposition}

\subsection{Focusing spaces and algebras}

The proofs of the properties stated in this section are given in
Costin (1997)\label{Costin++J}.

We say that a family of norms $\|\|_{\nu}$ depending on a parameter
$\nu\in\RR^+ $ is {\bf focusing} if for any $f$ with $\|f\|_{\nu_0}<\infty$

\begin{eqnarray}
  \label{focus-pocus}
  \|f\|_\nu\downarrow 0 \mbox{ as }\nu\uparrow\infty
\end{eqnarray}

Let $\mathcal{E}$ be a linear space and $\{\|\|_\nu\}$ a family of
norms satisfying (\ref{focus-pocus}).  For each $\nu$ we define a
Banach space $\mathcal{B}_\nu$ as the completion of
$\{f\in\mathcal{E}:\|f\|_{\nu}<\infty\}$. Enlarging $\mathcal{E}$ if
needed, we may assume that $\mathcal{B}_\nu\subset\mathcal{E}$. For
$\alpha<\beta$, (\ref{focus-pocus}) shows that the identity is an
embedding of $\mathcal{B}_\alpha$ in $\mathcal{B}_\beta$. Let
$\mathcal{F}\subset\mathcal{E}$ be the projective limit of the
$\mathcal{B}_\nu$.  That is to say

\begin{eqnarray}
  \label{focuproj}
  \mathcal{F}:=\bigcup_{\nu>0}\mathcal{B}_\nu
\end{eqnarray}

\z is endowed with the topology  in which a sequence is convergent if it
converges in {\em some} $\mathcal{B}_\nu$. We call $\mathcal{F}$
a {\bf focusing space}.

Consider now the case when
$\left(\mathcal{B}_{\nu},+,*,\|\|_\nu\right)$ are commutative Banach
algebras. Then $\mathcal{F}$ inherits a structure of a commutative
algebra, in which $*$ (``convolution'') is continuous. We say that
$\left(\mathcal{F},*,\|\|_\nu\right)$ is a {\bf focusing algebra}.

\subsection{Examples}
\label{sec:NC}

For more details see Costin (1997)\label{Costin++K}.  Let $K\in\RR^+$ and
$\mathcal{S}=\mathcal{S}_{K,\alpha_1,\alpha_2}=\{p:\arg(p)\in[\alpha_1,\alpha_2]\subset
(-\pi/2,\pi/2),|p|\le K\}$ (or a finite union of such sectors) and
$\mathcal{V}$ be a small neighborhood of the origin.
$\overline{\mathcal{V}}$ will be the closure of $\mathcal{V}$, cut
along the negative axis, and together with these upper and lower cuts.

 {\bf (1)}. $\lone_\nu(\mathcal{K})$. Let
 $\mathcal{K}=\mathcal{S}_{K,\phi,\phi}$.
 The space
$\lone_\nu(\mathcal{K})$ with the
 convolution $f*g:=p\mapsto\int_0^p f(s)g(p-s)\mathrm{d}s$  is a commutative
Banach algebra under each  of the (equivalent) norms

\begin{eqnarray}
  \label{norm00}
  \|f\|_\nu=\int_0^K \mathrm{e}^{-\nu t}|f(t\exp(\mathrm{i}\phi))|\mathrm{d}t
\end{eqnarray}

{\bf (2)} If $K=\infty$ in example $(1)$, then the norms (\ref{norm00})
are not equivalent anymore for different $\nu$, but convolution is
still continuous in (\ref{norm00}) and  the projective limit of
the $L^1_\nu(\RR^+ \mathrm{e}^{\mathrm{i}\phi})$, $\mathcal{F}(\RR^+
\mathrm{e}^{\mathrm{i}\phi})\subset L^1_{loc}(\RR^+ \mathrm{e}^{\mathrm{i}\phi})$, is a focusing algebra.

{\bf (3a)} $\mathcal{T}_\beta(\mathcal{S}\cup\overline{\mathcal{V}})$.
For $\Re(\beta)> 0$ and $\phi_1\ne\phi_2$,
this space is given by $\{f:f(p)=p^{\beta}F(p)\}$,
where $F$ is analytic in the interior
of $\mathcal{S}\cup\mathcal{V}$ and continuous in
its closure. We take the family of (equivalent) norms

\begin{eqnarray}
  \label{normF1}
  \|f\|_{\nu,\beta}=K\sup_{s\in
   \mathcal{S}\cup\overline{\mathcal{V}}}\left|\mathrm{e}^{-\nu p}f(p)\right|
\end{eqnarray}

\z It is clear that convergence of $f$ in $\|\|_{\nu,\beta}$ implies
uniform convergence of $F$ on compact sets in
$\mathcal{S}\cup\mathcal{V}$ (for $p$ near zero, this follows from
Cauchy's formula).  $\mathcal{T}_{\beta}$ are thus Banach spaces and
 focusing spaces by (\ref{normF1}). The spaces
$\{\mathcal{T}_{\beta}\}_{\beta}$ are isomorphic to each-other.
The application
\begin{eqnarray}
  \label{convodom1}
  (\cdot *\cdot):\mathcal{T}_{\beta_1}\times
\mathcal{T}_{\beta_2}\mapsto \mathcal{T}_{\beta_1+\beta_2+1}
\end{eqnarray}

\z is continuous:

A natural generalization of $\mathcal{T}_\beta$ is obtained taking
$\beta_1,\ldots,\beta_N\in\CC$ with positive real parts, no two of
them differing by an integer. If $f_\beta=\sum_{i=1}^k p^{\beta_i}A_i(p)$
with $A_i$ analytic, then $f_\beta\equiv 0$ iff $A_i\equiv 0$ for all $i$
(e.g., by a  Puiseux series argument).  It is then
natural to identify the space $\mathcal{T}_{\{\beta_1,\ldots,\beta_k\}}$
of functions of the form $f_\beta$ with
$\oplus_{i=1}^k \mathcal{T}_{\beta_i}$. Convolution with analytic
functions is defined on $\mathcal{T}_{\{\beta_1,\ldots,\beta_k\}}$ while
convolution of two functions in $\mathcal{T}_{\{\beta_1,\ldots,\beta_k\}}$
takes values in $\mathcal{T}_{\{\beta_i+\beta_j \,\mbox{mod}\ 1\} }$.
We write $\mathcal{T}_{\{\cdot\}}$ when the concrete
values of $\beta_1,\ldots,\beta_k$  do not matter.

{\bf (3b)} A particular case of the preceding example is
$\mathcal{A}_{z,l}(\mathcal{S\cup \mathcal{V}})$ consisting of
analytic functions in the interior of $\mathcal{S}\cup\mathcal{V}$,
continuous on its closure, and vanishing at the origin together with
the first $l$ derivatives. $\mathcal{A}_{z,l}$ can be identified with
$\mathcal{T}_l$.

{\bf (4)} $\mathcal{D}'_{m,\nu}$, the ``staircase distributions''.
 Let $\mathcal{D}(0,x)$ be the test
functions on $(0,x)$ and $\mathcal{D}=\mathcal{D}(0,\infty)$. Let
$\mathcal{D}'_m\subset\mathcal{D}'$ be the distributions $f$ for which
$f=F_k^{(km)}$ on $\mathcal{D}(0,k+1)$ with $F_k\in\lone[0,k+1]$.
There is a uniquely associated staircase decomposition, a sequence
$\left\{\Delta_i(f)\right\}_{i\in\NN}=\left\{\Delta_i\right\}_{i\in\NN}$
such that $\Delta_i\in\lone(\RR^+)$,
$\Delta_i=\Delta_i\bchi_{[i,i+1]}$ and

\begin{eqnarray}
  \label{stdec0}
  f=\sum_{i=0}^{\infty}\Delta_i^{(mi)}
\end{eqnarray}

\z With respect to the norm 

\begin{eqnarray}
  \label{normdistr00}
   \|f\|_{\nu ,m}:=\sqrt{2}\sum_{i=0}^{\infty}\nu ^{im}\|\Delta_i\|_\nu 
\end{eqnarray}

\z where $\|\Delta\|_\nu$ is computed from  (\ref{norm00}) with
$K=\infty$
and with convolution defined as

\begin{eqnarray}
  \label{formulaforconv}
   \label{stdeccomv}
  \Delta_k(f*\tilde{f})=\sum_{i+j=k}\Delta_i*\tilde{\Delta}_j-
\mathcal{P}^{m}\left\{\sum_{i+j=k-1}
\left(\Delta_i*\tilde{\Delta}_j\right)\bchi_{[0,k+1]}\right\}
\end{eqnarray}

\z $(\mathcal{D}'_m,+,*)$ is a commutative Banach algebra.
With respect to the family of norms $\|\|_{m,\nu}$, the projective
limit of the $\mathcal{D}'_{m,\nu}$, $\mathcal{F}_{m}$ is a
focusing algebra.

For any  $f\in L^1_{\nu_0}(\RR^+)$ there is a constant
$C(\nu,\nu_0)$ such that $f\in \mathcal{D}'_{m,\nu}$ for all
$\nu>\nu_0$ and

\begin{eqnarray}
  \label{majornorm}
  \|f\|_{\mathcal{D}'_{m,\nu}}\le C(\nu_0,\nu)\|f\|_{ L^1_{\nu_0}}
\end{eqnarray}

\z and formula (\ref{formulaforconv}) is equivalent to the usual convolution
in this case.

For $a\in\RR^+$,
$\mathcal{D}'_{m,\nu}(a,\infty)=\{f\in\mathcal{D}'_{m,\nu}:
\Delta_i(x)=0$ for  $x<a\}$ is a closed ideal in $\mathcal{D}'_{m,\nu}$
(isomorphic to the restriction $\mathcal{D}'_{m,\nu}(a,\infty)$ of
$\mathcal{D}'_{m,\nu}$ to $\mathcal{D}(a,\infty)$). The restrictions
$\mathcal{D}'_{m,\nu}(a,b)$ of $\mathcal{D}'_{m,\nu}$
to $\mathcal{D}(a,b)$ are for $0<a<b<\infty$ Banach spaces 
 with respect to the norm
(\ref{normdistr00}) restricted to $(a,b)$.

The functions in
$\mathcal{D}\left(\RR^+\backslash\NN\right)$ 
are dense in $\mathcal{D}'_{m,\nu}$, with respect to the norm
(\ref{normdistr00}).

If we choose a different interval length $l>0$ instead of $l=1$ in the
partition associated to (\ref{stdec0}), we then write
$\mathcal{D}'_{m,\nu}(l)$. Obviously, dilation gives a natural
isomorphism between these structures. If
$d=\{t\mathrm{e}^{\mathrm{i}\phi}:t\in\RR^+\}$ is any ray,
$\mathcal{D}'_{m,\nu}(d)$ and
$\mathcal{F}_{m;\phi}$ are defined in an analogous way and have
the same properties as their real counterpart.

Laplace transforms are naturally defined in
$\mathcal{D}'_{m,\nu}$.

\begin{Lemma}\label{existe} Laplace transform extends continuously
from $\mathcal{D}(\RR^+\backslash\NN)$ to $\mathcal{D}'_{m,\nu}(\RR^+)$ by the formula

\begin{eqnarray}
  \label{laptra}
  (\lap f)(x):=\sum_{k=0}^{\infty}x^{mk}\int_0^{\infty}\mathrm{e}^{-sx}\Delta_k(s)\mathrm{d}s
\end{eqnarray}

\z In particular, with $f,g, h'\in\mathcal{D}'_{m,\nu }$ we have
\begin{multline}
  \label{lapcomuta}
  \lap(f*g)=\lap(f)\lap(g)\cr 
 \lap(h')=x\lap(h)-h(0)\cr 
\lap(pf)=-(\lap(f))'
\end{multline}

For $x\in S_\nu= \{x:\Re(x)>\nu\} $ the sum (\ref{laptra}) converges
absolutely.  Laplace transform is,  for fixed $x\in S_\nu$, a continuous
functional (of norm less than one) on $\mathcal{D}'_{m,\nu}$.

$(\lap f)(x)$ is analytic in $S_\nu$.

Furthermore, $\mathcal{L}$  is {\em injective} in $\mathcal{D}'_{m,\nu}$. 
\end{Lemma}

\subsection{Appendix 2. Two examples of resonant equations}

Once resonant equations are allowed, Berry transitions
tend to become more complicated:
$y=\mbox{e$^{-x}$Ei}(x)+\exp({-x^2-\mbox{i}x^4})$ is a (least term
summable for $x\rightarrow +\infty$) solution of a linear
 differential equation with rational coefficients:

$$\left[\frac{d}{dx}-\frac{P'}{P}+2x+4\mbox{i}x^3\right]\left[x\frac{d^2}{dx^2}+(x+1)\frac{d}{dx}+1\right]y=0$$

\z where $P(x)=-16x^7+16\mbox{i}x^5-4\mbox{i}x^4 +(4-16\mbox{i})x^3
-2x^2-4x+1$. The term $\exp({-x^2-\mbox{i}x^4})$ is not seen by
least term truncation on the real axis but becomes much larger than
the error function contribution before entering the Berry region $\arg(x)\sim
|x|^{-1/2}$.

Borel summation using Ecalle acceleration (see \'Ecalle 1993) gives an
unambiguous description of the solutions of this type of equations and
can be used to decompose solutions conveniently before analyzing
transitions.  On the other hand, imposing that all Stokes constants
are nonzero suffices to exclude this and similar examples, but looks
rather restrictive and hard to verify.

 We next consider
the Berry transition of a family of \emph{resonant} equations with {\em
  nonzero} Stokes constants. The formal solutions of such equations
depart from Dingle's rule, and
also exhibit an interesting splitting of the Stokes rays, with two
Berry transitions in the same region. (We are mainly aiming at
illustration and the calculation is heuristic but a rigorous treatment
along these same lines is not difficult).  Let first $m\in\RR$ and

\begin{eqnarray}
  \label{e1}
  L[y]=y''+2y'+\left(1+\frac{m^2}{x}\right)y=\frac{1}{x}
\end{eqnarray}

\z {\sc (A)} {\em Formal solutions; behavior of coefficients.}

\smallskip

\smallskip

\z Taking $\tilde{y}=\sum_{k=0}^{\infty}a_kx^{-k}$ we get for $a_k$:

\begin{eqnarray}
  \label{r1}
 a_{k+1}=(2k-m^2)a_k-k(k-1)a_{k-1}
\end{eqnarray}

\z With $a_k=(k-1)!b_k$ we have

\begin{eqnarray}
  \label{r1b}
 b_{k+1}=(2-\frac{m^2}{k})b_k-b_{k-1}
\end{eqnarray}

\z which we analyze for $k\gg 1$ by WKB (an explicit solution is also
possible in this case). With $b_k=\mathrm{e}^{w_k}$ we get:

\begin{multline}
  \label{ewk}
  \e^{w+w'+\frac{1}{2}w''}\sim(2-\frac{m^2}{k})\e^w-\e^{w-w'+\frac{1}{2}w''}
\cr 
\Rightarrow\e^{\frac{1}{2}w''}\cosh(w')\sim 1-\frac{m^2}{2k}
(1+\frac{1}{2}w'')(1+\frac{1}{2}{w'}^2)\sim
1-\frac{m^2}{2k}
\cr\Rightarrow w'\sim\pm \i \sqrt{\frac{m^2}{k}+w''}
\sim\pm \i \sqrt{\frac{m^2}{k}\mp\i \frac{m}{2 k^{3/2}}}
\sim \pm \i \frac{m}{\sqrt{k}}+\frac{1}{4k}\end{multline}
\z so that 
\begin{multline}
b_k\sim k^{1/4}\left(A_+\e^{2\i m\sqrt{k}}+A_-\e^{-2\i
  m\sqrt{k}}\right)\mbox{ and } a_k=a_k^{+}+a_k^{-}
\cr \mbox{ with } a_{k}^{\pm}\sim (k-1)!A_{\pm} k^{1/4} e^{\pm 2\i m\sqrt{k}}
\end{multline}

\z and, in analogy with the nonresonant case we write
$\tilde{y}=\tilde{y}^++\tilde{y}^-$ where
$\tilde{y}^\pm=\sum_{k=0}^{\infty}a_k^{\pm}x^{-k}$.  We first see that
there is no curve $\arg(x)=f(|x|)$ for which the terms of
$\tilde{y}^{\pm}$ (or the terms of any linear combination of
$y^{\pm}$) have the same phase. There are instead two  parabolic
curves, $\pm\arg(x)=2m|x|^{-1/2}$ along which for $n\sim |x|$ we have
$[a_{n+1}^+/x^{n+1}]/[a_{n}^+/x^{n}] =1+o(1)$. This does not amount to
the terms being in phase, but rather is the discrete equivalent of a
stationary point.

 To reinterpret
Dingle's rule for this example, we will see
that there exist two Stokes parabolas, each associated
with one degree of freedom in the original equation
along which the transitions in the constants beyond all orders,
as measured by optimal truncation, are maximal.

\smallskip

\z {\sc (B)} {\em Asymptotic solutions of the homogeneous equation.}

\z Taking  $y=\e^w$ we obtain:

\begin{multline}
  \label{wk2}
  w''+{w'}^2+2w'+1+\frac{m^2}{x}=0\Rightarrow
  w'=-1\pm\i\sqrt{\frac{m^2}{x}
+w''}\cr w\sim -x\pm2\i m\sqrt{x}+\frac{1}{4}\ln x\Rightarrow
y_{\pm}\sim x^{1/4}\e^{-x\pm2\i m\sqrt{x}}
\end{multline}

\z as before.

\smallskip

\z {\sc (C)}{\em Berry smoothing.}

\z Let $y=\sum_{k=1}^{n-1}a_kx^{-k}+\sum_{\pm}C_{\pm}(x)y_\pm(x)$.
Then $L[y]-x^{-1}$ gives

\begin{multline}
  \label{apprx}
  \sum_{\pm}y_\pm\left[C_{\pm}''+2C'_\pm\left(1+\frac{y'_\pm}{y_\pm}\right)
\right]+\frac{n(n-1)a_{n-1}}{x^{n+1}}+\frac{(n-1)(n-2)a_{n-2}}{x^{n}}
\cr+\frac{m^2-2n-2}{x^n}a_{n-1}=L_+C_+ +L_- C_-+\frac{n(n-1)a_{n-1}}{x^{n+1}}
-\frac{a_n}{x^n}
\cr
=L_+C_+ +L_- C_-+\frac{n!}{x^n}b_{n}\left(\frac{1}{n}-\frac{b_{n-1}}{b_n x}\right)
\end{multline}

\z with obvious notations. Changing variables to $x=n\e^{\i\beta
  n^{-1/2}}$ above we get, to leading order,

\begin{multline}
  \label{eq3}
  n^{1/4}\sum_{\pm}\e^{-n+\i\sqrt{n}(\pm 2
    m-\beta)+\frac{1}{2}\beta^2\mp m\beta}\left(-\frac{C_{\pm}''}{n}\mp\frac{2m
C_{\pm}'}{n}\right)\cr=-\sum_{\pm}\frac{n!}{n^n}\e^{-\i\sqrt{n}\beta}
\left[n^{1/4}A_\pm\e^{\pm 2\i m n^{1/2}}n^{-3/2}\i (
\pm m+\beta)\right]\cr 
\end{multline}

\z Equating the coefficients of $\e^{2\i m n^{1/2}}$ we get,
with $B_{\pm}=Const+o(1)$
the system

\begin{eqnarray}
  \label{sy1}
  C_\pm''\pm\frac{2\i
  m}{x^{1/2}}C_\pm'=\i B_{\pm}n^{-1}\e^{-\frac{1}{2}\beta^2\pm m\beta}(\pm m+\beta)
\end{eqnarray}

\z where we change variables to $x=n\exp(\i \beta n^{-1/2})$
and get to leading order:

\begin{eqnarray}
  \label{sys2}
 -C_\pm''(\beta)\pm 2mC_\pm
'(\beta)=\i B_{\pm}\e^{-\frac{1}{2}\beta^2\pm m\beta}(\pm m+\beta)
\end{eqnarray}

\z with the bounded solutions:

\begin{eqnarray}
  \label{sol1}
  C_\pm(\beta)=\i \sqrt{\frac{\pi}{2}}\e^{\frac{1}{2}m^2}
\mbox{erf}\left(\frac{1}{\sqrt{2}}\beta-\frac{1}{\sqrt{2}}m\right)B_{\pm}+\mbox{Const.}_{\pm}
\end{eqnarray}

\begin{section}{Acknowledgments}

  One of the authors (OC) is very grateful to Professors Michael
  Berry, Percy Deift and Jean \'Ecalle for very interesting
  discussions.

\end{section}

\end{document}